\documentclass[12pt,a4paper]{article}
\usepackage{latexsym,epsfig,fullpage,amsfonts,amssymb,amsmath,amscd,epic,graphics,rotating,theorem,yhmath}
\usepackage[all]{xy}
\usepackage[usenames]{color}
\baselineskip 18pt \textwidth 15cm  \theoremstyle{plain}
\newtheorem{lemma}{Lemma}
\newtheorem{proposition}[lemma]{Proposition}\newtheorem{conjecture}{Conjecture}
\newtheorem{remark}[lemma]{Remark}

\newtheorem{theorem}{Theorem}

{\theorembodyfont{\rmfamily}
\begin{document}
\newcommand{\Pic}{{\operatorname{Pic}}}
\newcommand{\BH}{{\operatorname{bh}}}
\newcommand{\GW}{{\operatorname{GW}}}
\newcommand{\ini}{{\operatorname{ini}}}
\newcommand{\fm}{{\mathfrak{m}}}
\newcommand{\arc}{{\operatorname{arc}}}
\newcommand{\conv}{{\operatorname{Conv}}}
\newcommand{\Conj}{{\operatorname{conj}}}
\newcommand{\Ev}{{\operatorname{Ev}}}
\newcommand{\pr}{{\operatorname{pr}}}
\newcommand{\val}{{\operatorname{val}}}
\newcommand{\eps}{{\varepsilon}}
\newcommand{\idim}{{\operatorname{idim}}}
\newcommand{\DP}{{\operatorname{DP}}}
\newcommand{\Arc}{{\operatorname{Arc}}}
\newcommand{\Graph}{{\operatorname{Graph}}}
\newcommand{\sm}{{\mathrm{sm}}}
\newcommand{\Z}{{\mathbb Z}}
\newcommand{\bw}{{\boldsymbol{w}}}
\newcommand{\R}{{\mathbb{R}}}
\newcommand{\bz}{{\boldsymbol{z}}}
\newcommand{\bp}{{\boldsymbol{p}}}
\newcommand{\Hom}{{\operatorname{Hom}}}
\newcommand{\const}{{\operatorname{const}}}
\newcommand{\Tors}{{\operatorname{Tors}}}
\newcommand{\Spec}{{\operatorname{Spec}}}
\newcommand{\Aut}{{\operatorname{Aut}}}
\newcommand{\C}{{\mathbb{C}}}
\newcommand{\bA}{{\boldsymbol{A}}}\newcommand{\bF}{{\boldsymbol{F}}}
\newcommand{\bB}{{\boldsymbol{B}}}
\newcommand{\bn}{{\boldsymbol{n}}}
\newcommand{\bs}{{\boldsymbol{s}}}
\newcommand{\PP}{{\mathbb{P}}}
\newcommand{\Id}{{\operatorname{Id}}}
\newcommand{\Sing}{{\operatorname{Sing}}}
\newcommand{\codim}{{\operatorname{codim}}}
\newcommand{\Ann}{{\operatorname{Ann}}}
\newcommand{\ord}{{\operatorname{ord}}}
\newcommand{\mt}{{\operatorname{mt}}}\newcommand{\br}{{\operatorname{br}}}
\newcommand{\Prec}{{\operatorname{Prec}}}
\newcommand{\proofend}{\hfill$\blacksquare$\bigskip}
\newcommand{\pperp}{{\perp{\hskip-0.4cm}\perp}}
\newcommand{\bk}{{\boldsymbol{k}}}
\newcommand{\bl}{{\boldsymbol{l}}}

\title{On Welschinger invariants of descendant type}


\author{Eugenii Shustin}

%
%
%
\date{}
\maketitle
\bigskip

\centerline{\it Dedicated to Gert-Martin Greuel in occasion of his 70-th birthday}

\medskip

\begin{abstract}
We introduce enumerative invariants
of real del Pezzo surfaces that count real rational curves
belonging to a given divisor class, passing through a generic
conjugation-invariant configuration of points and satisfying
preassigned tangency conditions to given smooth arcs centered at the fixed points.
The counted curves are equipped with Welschinger-type signs. We prove
that such a count does not depend neither on the choice of the point-arc
configuration, nor on the variation of
the ambient real surface. These invariants can be regarded as
a real counterpart
of (complex) descendant invariants.
\end{abstract}

\maketitle

\section*{Introduction}

Welschinger invariants of real rational symplectic manifolds
\cite{Welschinger:2003,Welschinger:2005,Welschinger:2005a,
Welschinger:2009}
serve as genus zero open Gromov-Witten invariants.
In dimension four and in the algebraic-geometric setting, they are
well-defined for real del Pezzo surfaces
(cf. \cite{Itenberg_Kharlamov_Shustin:2012}), and they count real
rational curves in a given divisor class
passing through a generic conjugation-invariant configuration of
points, and equipped with weights $\pm1$. An important outcome of
Welschinger's theory is that, whenever
Welschinger invariant does not vanish, there exists a real rational
curve of a given divisor class matching
an appropriate number of arbitrary generic conjugation-invariant constraints.

There are several extensions of the original Welschinger invariants:
modifications for multi-component
real del Pezzo surfaces
\cite{Itenberg_Kharlamov_Shustin:2013a,Itenberg_Kharlamov_Shustin:2012},
mixed and relative invariants \cite{IKS,RS,We5}, invariants of positive
genus for multi-component
real del Pezzo surfaces \cite{Sh4} and for ${\mathbb P}^{2k+1}$, $k\ge1$) \cite{GZ1,GZ2}.
The goal of this paper is to introduce
Welschinger-type invariants for
real del Pezzo surfaces, which count real rational curves passing through
some fixed points and tangent
to fixed smooth arcs centered at the fixed points. They can be viewed as a real
counterpart of certain descendant invariants (cf. \cite{GKP}).

The main result of this note is Theorem \ref{t1} in Section \ref{sec1},
which states the existence of invariants
independent of the choice of constraints and of the variation of the surface.
Our approach in general is similar to that in \cite{Itenberg_Kharlamov_Shustin:2012}, and
it consists in the study of codimension one bifurcations of the set of curves subject to
imposed constraints when one varies either the constraints, or the real and complex
structure of the surface. In Section \ref{sec2}, we show a few simple examples.
The computational aspect and quantitative properties of the invariants will be treated in a forthcoming paper.

\section{Invariants}\label{sec1}

Let $X$ be a real del Pezzo surface with a nonempty
real point set $\R X$, and $F\subset\R X$ a connected component.
Pick a conjugation-invariant class $\varphi\in H_2(X\setminus F;\Z/2)$.
Denote by $\Pic^\R_+(X)\subset\Pic(X)$ the
subgroup of real effective divisor classes. Pick a non-zero class $D\in\Pic^\R(X)$,
which is $F$-compatible
in the sense of \cite[Section 5.3]{Itenberg_Kharlamov_Shustin:2013}. Observe that any real
rational (irreducible) curve
$C\in|D|$ has a one-dimensional real branch (see, for example,
\cite[Section 1.2]{Itenberg_Kharlamov_Shustin:2012}), and hence we can define $C_+,C_-$,
the images of the components of $\PP^1\setminus\R\PP^1$ by the normalization map.

Given a smooth (complex) algebraic variety $\Sigma$, a point $z\in
\Sigma$, and a positive integer $s$, the space of $s$-arcs in $\Sigma$ at $z$ is
$$\Arc_s(\Sigma,z)=\Hom(\Spec\C[t]/(t^{s+1}),(\Sigma,z))/\Aut(\C[t]/(t^{s+1}))\ .$$
Denote by $\Arc_s^{\sm}(\Sigma,z)\subset\Arc_s(\Sigma,z)$ the (open) subset
consisting of smooth
$s$-arcs, i.e., of those which are represented by an embedding $(\C,0)\to(\Sigma,z)$.

Choose two collections of positive integers
$\bk=\{k_i,\ 1\le i\le r\}$ and $\bl=\{l_j,\ 1\le j\le m\}$, where $r,m\ge0$ and
\begin{equation}\sum_{i=1}^rk_i+2\sum_{j=1}^ml_j=-DK_X-1\ ,\label{e2}\end{equation}
and all $k_1,...,k_r$ are odd. Pick distinct points
$z_1,...,z_r\in F$ and real arcs $\alpha_i\in\Arc_{k_i}^{\sm}(X,z_i)$, $1\le i\le r$,
and also distinct points
$w_1,...,w_m\in X\setminus\R X$ and arcs $\beta_j\in\Arc_{l_j}^{\sm}(X,w_j)$. Denote
$\bz=(z_1,...,z_r)$,
$\bw=(w_1,\overline w_1,...,w_m,\overline w_m)$ and
\begin{equation}{\mathcal A}=(\alpha_1,...,\alpha_r)\in\prod_{i=1}^r
\Arc_{k_i}^{\sm}(X,z_i)\ ,\label{bA}\end{equation}
\begin{equation}
{\mathcal B}=(\beta_1,\overline\beta_1,...,\beta_m,\overline\beta_m)\in
\prod_{j=1}^m\left(\Arc_{l_j}^{\sm}(X,w_j)\times\Arc_{l_j}^{\sm}(X,\overline w_j)\right)\ .
\label{bB}\end{equation} In the moduli space ${\mathcal M}_{0,r+2m}(X,D)$ of stable maps of
rational curves with $r+2m$ marked points, we consider the subset
${\mathcal M}_{0,r+2m}(X,D,(\bk,\bl),(\bz,\bw),({\mathcal A},
{\mathcal B}))$ consisting of the
elements $[\bn:\PP^1\to X,\bp]$, $\bp=(p_1,...,p_r,q_1,...,q_m,
q'_1,...,q'_m)\subset\PP^1$, such that
$$\bn^*\left(\bigcup{\mathcal A}\cup\bigcup{\mathcal B}\right)\ge
\sum_{i=1}^rk_ip_i+\sum_{j=1}^ml_j(q_j+q'_j)\ .$$
Let ${\mathcal M}_{0,r+2m}^{\;im,\R}(X,D,(\bk,\bl),(\bz,\bw),({\mathcal A},
{\mathcal B}))\subset{\mathcal M}_{0,r+2m}(X,D,(\bk,\bl),(\bz,\bw),({\mathcal A},
{\mathcal B}))$ be the set of elements $[\bn:\PP^1\to
X,\bp]$ such that $\bn$ is a conjugation invariant immersion, the points
$p_1,...,p_r\in\PP^1$ are real,
and $q_j,q'_j\in\PP^1$ are complex conjugate, $j=1,...,m$. For a generic
choice of
point sequences $\bz$ and $\bw$, and arc sequences
${\mathcal A}$ and ${\mathcal B}$ in the arc spaces indicated in (\ref{bA})
and (\ref{bB}), the set
${\mathcal M}_{0,r+2m}^{\;im,\R}(X,D,(\bk,\bl),(\bz,\bw),({\mathcal A},
{\mathcal B}))$ is finite (cf. Proposition \ref{p4}(1) below).

Given an element $\xi=[\bn:\PP^1\to X,\bp]\in{\mathcal M}_{0,r+2m}^{\;im,\R}
(X,D,(\bk,\bl),(\bz,\bw),({\mathcal A},
{\mathcal B}))$, denote $C=\bn(\PP^1)$ and
define the {\it Welschinger sign} of $\xi$ by (cf. \cite[Formula (1)]{Itenberg_Kharlamov_Shustin:2012})
$$W_\varphi(\xi)=(-1)^{C_+\circ\; C_-+C_+\circ\;\varphi}\ .$$ Notice that, if
$C$ is nodal, then $C_+\circ C_-$ has the same parity as the number of real
solitary nodes of $C$ (i.e. nodes locally equivalent to $x^2+y^2=0$).

Finally, put
\begin{equation}W(X,D,F,\varphi,(\bk,\bl),(\bz,\bw),({\mathcal A},
{\mathcal B}))=\sum_{\xi\in{\mathcal M}_{0,r+2m}^{\;im,\R}(X,D,(\bk,\bl),
(\bz,\bw),({\mathcal A},{\mathcal B}))}W_\varphi(\xi)\ .\label{e1}\end{equation}

\begin{theorem}\label{t1}
(1) Let $X$ be a real del Pezzo surface with $\R X\ne\emptyset$, $F\subset\R X$ a
connected component,
$\varphi\in H_2(X\setminus F,\Z/2)$ a conjugation-invariant class, $D\in\Pic^\R_+(X)$
a nef and big,
$F$-compatible divisor class, 
$\bk=(k_1,...,k_r)$ a (possibly empty) sequence of
positive odd integers such that
\begin{equation}\max\{k_1,...,k_r\}\le3\ ,\label{e33}\end{equation} and
$\bl=(l_1,...,l_m)$ a (possibly empty) sequence of positive integers
satisfying (\ref{e2}), $\bz=(z_1,...,z_r)$ a sequence of distinct points of $F$, $\bw=(w_1,...,w_m,
\overline w_1,...,\overline w_m)$ a sequence of distinct points of $X\setminus\R X$, and,
at last, ${\mathcal A}$, ${\mathcal B}$ are
arc sequences as in (\ref{bA}), (\ref{bB}). Then the number
$W(X,D,F,\varphi,(\bk,\bl),(\bz,\bw),({\mathcal A},
{\mathcal B}))$ does not depend neither on the choice of generic point configuration
$\bz$, $\bw$, nor on the choice of arc sequences ${\mathcal A}$, ${\mathcal B}$ subject to
conditions indicated above.

(2) If tuples $(X,D,F,\varphi)$ and $(X',D',
F',\varphi')$ are deformation equivalent so that $X$ and $X'$ are joined by
a flat family of real smooth rational surfaces,
then we have (omitting $(\bz,\bw)$ and $({\mathcal A},{\mathcal B})$ in the notation)
$$W(X,D,F,\varphi,(\bk,\bl))=W(X',D',F',\varphi',(\bk,\bl))
\ .
$$
\end{theorem}

\begin{remark}\label{r1} (1) If $k_i=l_j=1$ for all $1\le i\le r$, $1\le j\le m$, then we obtain
original Welschinger invariants in their modified form \cite{Itenberg_Kharlamov_Shustin:2013a},
and hence the required statement follows from
\cite[Proposition 4 and Theorem 6]{Itenberg_Kharlamov_Shustin:2012}.
This, in particular, yields that we have to consider the only case $-DK_X-1\ge3$.

(2) In general, one cannot impose even tangency conditions at real points $z_1,...,z_r$. Indeed,
suppose that $r\ge1$ and $k_1=2s$ is even. Suppose that $-DK_X-1\ge2s$ and $p_a(D)=(D^2+DK_X)/2+1\ge s$.
In the linear system $|D|$, the curves, which
intersect the arc $A_1$ at $z_1$ with multiplicity $\ge s$ and have at least $s$ nodes,
form a subfamily of codimension $3s$.
On the other hand, the family of curves having singularity $A_{2s}$ at $z_1$ and $(s-1)$ additional
infinitely near to $z_1$ points lying on the arc $\alpha_1$, has codimension $3s+1$ and it
lies in the boundary of
the former family. Over the reals, this wall-crossing results in the change of the
Welschinger sign of the curve that undergoes the corresponding bifurcation. Indeed,
take local coordinates $x,y$ such that
$z_1=(0,0)$ and $\alpha_1=\{y=0\}$, and consider the family of curves
$$y=t^{2s},\quad x=\eps t+t^2+t^3,\quad \eps\in(\R,0)\ .$$ For $\eps=0$, the curve has singularity
$A_{2s}$ at $z_1$, and its next $(s-1)$ infinitely near to $z_1$ points belong to $\alpha_1$.
In turn, for $\eps\ne0$, the node, corresponding to the values $t=\pm\sqrt{-\eps}$, is solitary
as $\eps>0$ and non-solitary as $\eps<0$, whereas the remaining $(s-1)$ nodes stay imaginary or solitary.
\end{remark}

\begin{conjecture}\label{conj1}
Theorem \ref{t1} is valid without restriction (\ref{e33}).
\end{conjecture}

The proof of Theorem \ref{t1} in general follows the lines of \cite{Itenberg_Kharlamov_Shustin:2012},
where we verify the constancy of the introduced enumerative numbers
in one-dimensional families of constraints and families of surfaces. The former verification requires
a classification of codimension one degenerations of the curves in count, while the
latter verification is based on a suitable analogue of the Abramovich-Bertram-Vakil
formula \cite{Abramovich_Bertram:2001,Vakil:2000}. Restriction (\ref{e33}) results from the lack
of our understanding of non-reduced degenerations of the counted curves.

\section{Degeneration and deformation of curves on complex rational surfaces}\label{sec-aa}

\subsection{Auxiliary miscellanies}\label{sec2.1}
{\bf(1) Tropical limit.} For the reader's convenience, we shortly remind what is the tropical limit
in the sense of \cite[Section 2.3]{Sh0}, which will be used below. In the field of
complex Puiseux series $\C\{\{t\}\}$, we consider the non-Archimedean valuation
$\val(\sum_ac_at^a)=-\min\{a\ :\ c_a\ne0\}$. Given a polynomial (or a power series)
$F(x,y)=\sum_{(i,j)\in\Delta}c_{ij}x^iy^j$ over $\C\{\{t\}\}$ with Newton polygon $\Delta$,
its tropical limit consists of the following data:
\begin{itemize}\item a convex piece-wise linear function $N_F:\Delta\to\R$, whose graph
is the lower part of the polytope
$\conv\{(i,j,-\val(c_{ij}))\ :\ (i,j)\in\Delta\}$, the subdivision $S_F$ of $\Delta$ into
linearity domains of $N_F$, and the tropical curve $T_F$, the closure of
$\val(F=0)$;
\item limit polynomials (power series) $F^\delta_{\ini}(x,y)=\sum_{(i,j)\in\delta}
c_{ij}^0x^iy^j$, defined for any face $\delta$ of the subdivision $S_F$, where
$c_{ij}=t^{N_F(i,j)}(c_{ij}^0+O(t^{>0}))$ for all $(i,j)\in\Delta$.
\end{itemize}

\smallskip

\noindent{\bf(2) Rational curves with Newton triangles.}

\begin{lemma}\label{le7}
(1) For any integer $k\ge1$, there are exactly $k$ polynomials $F(x,y)=\sum_{i,j}c_{ij}x^iy^j$ with Newton
triangle $T=\conv\{(0,0),(0,2),(k,1)\}$, whose coefficients $c_{00},c_{01},c_{02},c_{11}$ are given
generic non-zero constants
and which define plane rational curves. Furthermore, in the space of polynomials with Newton triangle
$T$, the family
of polynomials defining rational curves intersects transversally with the linear subspace given by
assigning generic nonzero constant values to the coefficients $c_{00},c_{01},c_{02},c_{11}$.
If the coefficients $c_{00},c_{01},c_{02},c_{11}$ are real, then,
\begin{itemize}\item
for an odd $k$, there is an odd number of real
polynomials $F$ defining rational curves, and each of these curves has an even number of real solitary
nodes,
\item for an even $k$ there exists an even number (possibly
zero) of polynomials $F$ defining rational curves,
and half of these curves have an odd number of real solitary nodes while the other half an even number of
real solitary nodes.
\end{itemize}

(2) For any integer $k\ge1$, there are exactly $k$ polynomials $F(x,y)=\sum_{i,j}c_{ij}x^iy^j$ with Newton
triangle $T=\conv\{(0,0),(0,2),(k,1)\}$, whose coefficients $c_{00},c_{02},c_{11}$ are given
generic non-zero constants, the coefficient $c_{k-1,1}$ vanishes, and which define plane rational curves. Furthermore,
in the space of polynomials with Newton triangle
$T$ and vanishing coefficient $c_{k-1,1}$, the family
of polynomials defining rational curves intersects transversally with the linear subspace given by
assigning generic nonzero constant values to the coefficients $c_{00},c_{02},c_{11}$.
If the coefficients $c_{00},c_{02},c_{11}$ are real, then,
\begin{itemize}\item
for an odd $k$, there is a unique real
polynomial $F$ defining a rational curve, and this curve either has $k-1$ real solitary nodes,
or has no real nodes at all,
\item for an even $k$, either there are no real polynomials defining rational curves, or there are
two real polynomials, one defining a rational curve with $k-1$ real solitary nodes, and the other defining a rational curve without
real solitary nodes.
\end{itemize}
\end{lemma}

{\bf Proof.}
Both statements can easily be derived from \cite[Lemma 3.9]{Sh0}.
\proofend

\smallskip

\noindent
{\bf(3) Deformations of singular curve germs.} Our key tool in
the estimation of dimension of families of curves
will be \cite[Theorem 2]{GuS}
(see also \cite[Lemma II.2.18]{GLS}).
For the reader's convenience, we remind it here. Let $C$ be a reduced curve on a smooth
surface $\Sigma$, and $z\in C$. By $\mt(C,z)$ we denote we denote the intersection multiplicity at
$z$ of $C$ with a generic
smooth curve on $\Sigma$ passing through $z$, by $\delta(C,z)$ the $\delta$-invariant,
and by $\br(C,z)$ the number of irreducible components of $(C,z)$.

\begin{lemma}\label{gus}
Let $C_t$, $t\in(\C,0)$, be a flat family of reduced curves on a smooth surface $\Sigma$, and
$z_t\in C_t$, $t\in(\C,0)$, a section such that the family of germs $(C_t,z_t)$, $t\in(\C,0)$,
is equisingular. Denote by $U$ a neighborhood of $z_0$ in $\Sigma$, and by $(C\cdot C')_U$ the
total intersection number of curves
$C,C'$ in $U$. The following lower bounds hold:
\begin{enumerate}\item[(i)] $(C_0\cdot C_t)_U\ge\mt(C_0,z_0)-\br(C_0,z_0)+2\delta(C_0,z_0)$ for $t\in(\C,0)$;
\item[(ii)] If $L$ is a smooth curve passing through $z_0=z_t$, $t\in(\C,0)$, and $(C_t\cdot L)_{z_0}
=\const$,
then $$(C_0\cdot C_t)_U\ge(C_0\cdot L)_{z_0}+\mt(C_0,z_0)-\br(C_0,z_0)+2\delta(C_0,z_0)$$ for $t\in(C,0)$.
\item[(iii)] If $L$ is a smooth curve containing the family $z_t$, $t\in(\C,0)$, and $(C_t\cdot L)_{z_t}
=\const$,
then $$(C_0\cdot C_t)_U\ge(C_0\cdot L)_{z_0}-\br(C_0,z_0)+2\delta(C_0,z_0)$$ for $t\in(C,0)$.
\end{enumerate}
\end{lemma}

Let $x,y\in(\C,0)$ be local coordinates in a
neighborhood of
a point $z$ in a smooth projective surface $\Sigma$. Let $L=\{y=0\}$, and $(C,z)\subset(\Sigma,z)$
a reduced, irreducible curve germ such that
$(C\cdot L)_z=s\ge1$. Denote by $\fm_z\subset{\mathcal O}_{\Sigma,z}$ the maximal ideal and
introduce the ideal
$I^{L,s}_{\Sigma,z}=\{g\in\fm_z\ :\ \ord g\big|_{L,z}\ge s\}$. The semiuniversal
deformation base of the germ $(C,z)$ in the space of germs $(C',z)$ subject to
condition $(C'\cdot L)_z\ge s$ can be identified with the germ at zero of the space
$$B_{C,z}(L,s):=I^{L,s}_{\Sigma,z}\big/\big\langle f,\frac{\partial f}{\partial x}\cdot\fm_z,
\frac{\partial f}{\partial y}\cdot I_{\Sigma,z}^{L,s}\big\rangle\ ,$$
where $f\in{\mathcal O}_{\Sigma,z}$ locally defined the germ $(C,z)$
(cf. \cite[Corollary II.1.17]{GLS}).

\begin{lemma}\label{def}
(1) The stratum $B^{eg}_{C,z}(L,s)\subset B_{C,z}(L,s)$ parameterizing equigeneric deformations
of $(C,z)$ is smooth of codimension $\delta(C,z)$ and its tangent space is
\begin{equation}T_0B^{eg}_{C,z}(L,s)=I^{L,s}_{C,z}\big/\big\langle f,\frac{\partial f}{\partial x}\cdot\fm_z,
\frac{\partial f}{\partial y}\cdot I^{L,s}_{\Sigma,z}\big\rangle\ ,\label{edef1}\end{equation}
where
$$I^{L,s}_{C,z}=\{g\in{\mathcal O}_{\Sigma,z}\ :\ \ord g\big|_{C,z}\ge s+2\delta(C,z)\}\ .$$

(2) If $\Sigma$, $(C,z)$, and $L$ are real, and $s$ is odd, then a generic member of
$B^{eg}_{C,z}(L,s)$ is smooth at $z$ and has only imaginary and real solitary nodes;
the number of solitary nodes is $\delta(C,z)\mod2$.
\end{lemma}

{\bf Proof.} (1) In \cite[Lemma 13]{IKS}, we proved a similar statement for the case $s=2$ and
$(C,z)$ of type $A_{2k}$, $k\ge1$, and we worked with equations. Here we settle the general case,
and we work with parameterizations. First, observe that a general member of
$B^{eg}_{C,z}(L,s)$ has $\delta(C,z)$ nodes as its singularities and is smooth at
$z$. Thus, $\codim_{I^{L,s}_{\Sigma,z}}B^{eg}_{C,z}(L,s)=\delta(C,z)$, the tangent space to
$B^{eg}_{C,z}(L,s)$ at its generic point $C'$ is formed by the elements $g\in{\mathcal O}_{\Sigma,z}$,
which vanish at the nodes of $C'$ and whose restriction to $(L,z)$ has order $s$. Clearly, the
limits of these tangent
spaces as $C'\to(C,z)$ contain the space $I^{L,s}_{C,z}\big/\langle f,\frac{\partial f}{\partial x}\fm_z,
\frac{\partial f}{\partial y}I^{L,s}_{\Sigma,z}\rangle$. On the other hand,
$\dim I^{L,s}_{\Sigma,z}/I^{L,s}_{C,z}=\delta(C,z)$ (see, for example, \cite[Lemma 6]{Sh}).
Let us
show the smoothness of $B^{eg}_{C,z}(L,s)$. Notice that the germ $(C,z)$ admits a uniquely defined
parametrization $x=t^s$, $y=\varphi(t)$, $t\in(\C,0)$, where $\varphi(0)=0$, and each
element $C'\in B^{eg}_{C,z}(L,s)$ admits a unique parametrization
$x=t^s$, $y=\varphi(t)+\sum_{i=1}^ma_it^i$, where $m=\dim B^{eg}_{C,z}(L,s)$, $a_1,...,a_m\in(\C,0)$.
Choose $m$ distinct generic values $t_1,...,t_m\in(\C,0)\setminus\{0\}$ and take the germs of the lines
$L_i=\{(t_i^s,y)\ :\ y\in(\C,\varphi(t_i)\}$, $i=1,...,m$.
It follows that the stratum $B^{eg}_{C,z}(L,s)$ is diffeomorphic to $\prod_{i=1}^mL_i$;
hence the smoothness and (\ref{edef1}).

(2) The second statement follows from the observation that the equation
$t_1^s=t_2^s$ has no real solutions $t_1\ne t_2$.
\proofend

Let $C^{(1)},C^{(2)}\subset
\Sigma$ be
two distinct immersed rational curves, $z\in C^{(1)}
\cap C^{(2)}$
a smooth point of both $C^{(1)}$ and $C^{(2)}$, and $W_z\subset\Sigma$ a sufficiently
small neighborhood of $z$.
Denote by $V\subset|C^{(1)}+C^{(2)}|$ the germ at
$C^{(1)}\cup C^{(2)}$ of the family of curves, whose total $\delta$-invariant in $\Sigma\setminus U$
coincides with that of $C^{(1)}\cup C^{(2)}$.

\begin{lemma}\label{le4}
(1) The germ $V$ is smooth of dimension $$c=(C^{(1)}\cdot C^{(2)})_z-C^{(1)}K_\Sigma-C^{(2)}K_\Sigma-2
\ ,$$ and
its tangent space
isomorphically projects onto the space
${\mathcal O}_{\Sigma,z}/I_z$, where $$I_z=\{f\in{\mathcal O}_{\Sigma,z}\ :
\ \ord f\big|_{(C^{(i)},z)}\ge(C^{(1)}\cdot C^{(2)})_z-C^{(i)}K_\Sigma-1,\ i=1,2\}
\ .$$

(2) Let $f_1,...,f_c,f_{c+1},...$ be a basis of the tangent space to
$|C^{(1)}+C^{(2)}|$ at $C^{(1)}\cup C^{(2)}$ such that $f_1,...,f_c$ project to a basis of
${\mathcal O}_{\Sigma,z}/I_z$, and $f_j\in I_z$, $j>c$, satisfy
$$\ord f_{c+1}\big|_{(C^{(1)},z)}=(C^{(1)}\cdot C^{(2)})_z-C^{(1)}K_\Sigma-1\ ,$$
$$\ord f_j\big|_{(C^{(1)},z)}\ge(C^{(1)}\cdot C^{(2)})_z-C^{(1)}K_\Sigma,\quad
j>c+1\ ,$$
and let
$$\sum_{i=1}^ct_if_i+\sum_{j>c}a_j(\overline t)f_j, \quad \overline t=(t_1,...,t_c)\in(\C^{\;c},0)\ ,$$
be a parametrization of $V$, where $C^{(1)}\cup C^{(2)}$ corresponds to the origin, and
$a_j$, $j>c$, are analytic functions vanishing at zero. Then
\begin{equation}\frac{\partial a_{c+1}}{\partial t_i}(0)\ne0\quad\text{for all}\ 1\le i\le c\
\text{with}\ \ord f_i\big|_{(C^{(1)},z)}\le(C^{(1)}\cdot C^{(2)})_z
-C^{(1)}K_\Sigma-2\ .
\label{ered13}\end{equation}
\end{lemma}

{\bf Proof.} Let $\nu^{(i)}:\PP^1\to C^{(i)}\hookrightarrow\Sigma$ be the normalization,
$p_i=(\nu^{(i)})^*(z)$, $i=1,2$.
Note that by Riemann-Roch
$$h^k(\PP^1,{\mathcal N}^{\nu^{(i)}}_{\PP^1}(-(-C^{(i)}K_\Sigma-1)p_i))=0,\quad k=0,1,\ i=1,2\ ,$$
where
${\mathcal N}$ denotes the normal bundle of the corresponding map,
and observe that the codimension of $I_z$ in ${\mathcal O}_{\Sigma,z}$ equals $c$.
The first statement of lemma follows.

For the second statement, we note that a generic irreducible element $C\in V$ satisfies
\begin{eqnarray}(C\cdot C^{(1)})_{W_z}&\le &C^{(1)}C^{(2)}+(C^{(1)})^2-(C^{(1)}C^{(2)}-(C^{(1)}\cdot
C^{(2)})_z)\nonumber\\
& & -((C^{(1)})^2+C^{(1)}K_\Sigma+2)=(C^{(1)}\cdot C^{(2)})_z-C^{(1)}K_\Sigma-2\ .
\label{ered14}\end{eqnarray} Next, we choose $i\in\{1,...,c\}$ as in (\ref{ered13}) and
take $C\in V$ given by
the parameter values $t_i=t$, $t_j=t^s$ with some $s>1$ for all $j\in\{1,...,c\}\setminus\{i\}$.
Then,
if $a_{c+1}=O(t^m)$ with $m>1$, one encounters at least $(C^{(1)}\cdot C^{(2)})_z-C^{(1)}K_\Sigma-1$
intersection points
of $C$ and $C^{(1)}$ in $W_z$. Thus, (\ref{ered13}) follows.
\proofend

\smallskip

\noindent
{\bf (4) Geometry of arc spaces.} Let $\Sigma$ be a smooth projective surface. Given an
integer $s\ge0$, denote by
$\Arc_s(\Sigma)$ the vector bundle of $s$-arcs over $\Sigma$
and by $\Arc_s^{\sm}(\Sigma)$ the bundle of smooth $s$-arcs over $\Sigma$ (particularly,
$\Arc_0(\Sigma)=\Arc_0^{\sm}(\Sigma)=\Sigma$).
For any smooth curve $C\subset\Sigma$, we have a natural map $\arc_s:C\to\Arc_s^{\sm}(\Sigma)$, sending
a point $z\in C$ to the $s$-arc at $z$ defined by the germ $(C,z)$. The following statement
immediately follows from
basic properties of ordinary analytic differential equations:

\begin{lemma}\label{arc}
Let $s\ge1$, $U$ a neighborhood of a point $z\in\Sigma$, and $\sigma$ a smooth section of the
natural projection
$\pr_s:\Arc_s^{\sm}(U)\to\Arc_{s-1}^{\sm}(U)$. Then there exists  a smooth analytic curve
$\Lambda$ passing through $z$, defined in a neighborhood
$V\subset U$ of $z$, and such that $\arc_s(\Lambda)\subset\sigma(\Arc_{s-1}^{\sm}(V))$.
\end{lemma}

Now, let $\Sigma$ be a smooth rational surface, $\bn:\PP^1\to\Sigma$ an immersion,
$C=\bn(\PP^1)\in|D|$, where $-DK_\Sigma=k>0$.
Pick a point $p\in\PP^1$ such that $z=\bn(p)$ is a smooth point
of $C$.
Denote by $U\subset\Arc_{k-1}(\Sigma)$
the natural image of the germ of ${\mathcal M}_{0,1}(\Sigma,D)$ at $[\bn:\PP^1\to\Sigma,p]$.
Choose coordinates $x,y$ in a neighborhood of $z$ so that
$z=(0,0)$, $C=\{y+x^k=0\}$, and introduce two one-parameter subfamilies
$\Lambda'=(z'_t,\alpha'_t)_{t\in(\C,0)}$ and $\Lambda''=
(z''_t,\alpha''_t)_{t\in(\C,0)}$ of $\Arc_{k-1}(\Sigma)$:
$$z'_t=(t,0),\ \alpha'_t=\{y=(x-t)^l\},\quad z''_t=(0,0),\ \alpha''_t=\{y=tx^{k-1}\},
\quad t\in(\C,0)\ ,$$ where $l>k$.

\begin{lemma}\label{le6}
The germ $U$ is smooth of codimension one in $\Arc_{k-1}(\Sigma)$, and it transversally intersects both
$\Lambda'$ and $\Lambda''$.
\end{lemma}

{\bf Proof.}
It follows from Riemann-Roch and from Lemma
\ref{gus}(iii) that $V$ admits the following parametrization:
$$((x_0,y_0),\{y=y_0+a_1(x-x_0)+...+a_{k-2}(x-x_0)^{k-2}+\varphi(x_0,y_0,
\overline a)(x-x_0)^{k-1}\})\ ,$$
$$x_0,y_0,a_1,...,a_{k-2}\in(\C,0),\quad
\overline a=(a_1,...,a_{k-2}),\quad \varphi(0)=0,
\quad\frac{\partial\varphi}{\partial x_0}(0)\ne0\ .$$
Thus, $V$ is a smooth hypersurface. The required intersection transversality
results from a routine computation.
\proofend

\subsection{Families of curves and arcs on arbitrary del Pezzo surfaces}\label{sec2.2}
Let $\Sigma$ be a smooth del Pezzo surface of degree $1$, and
$D\in\Pic(\Sigma)$ be an effective divisor such that $-DK_\Sigma-1>0$. Fix positive integers
$n\le-DK_\Sigma-1$ and \mbox{$s\gg-DK_\Sigma-1$}.
Denote by $\ring\Sigma^n\subset\Sigma^n$ the complement of the diagonals and by
$\Arc_s(\ring\Sigma^n)$ the total space of the restriction to $\ring\Sigma^n$ of the bundle
$(\Arc_s(\Sigma))^n\to\Sigma^n$.
In this section, we stratify the space
$\Arc_s(\ring\Sigma^n)$ with respect to the intersection of arcs with rational curves
 in $|D|$, and we describe all strata of codimension zero and one.

Introduce the following spaces of curves: given $(\bz,{\mathcal A})\in\Arc_s(\ring\Sigma^n)$,
$\bz=(z_1,...,z_n)$, ${\mathcal A}=(\alpha_1,...,\alpha_n)$,
and a sequence $\bs=(s_1,...,s_n)\in\Z_{>0}^n$ summing up to $|\bs|\le s$, put
\begin{eqnarray}{\mathcal M}_{0,n}(\Sigma,D,\bs,\bz,{\mathcal A})&=&\{[\bn:\PP^1\to
\Sigma,\bp]\in{\mathcal M}_{0,n}(\Sigma,D)\ :\ \nonumber\\
& &\quad\bn(p_i)=z_i,\quad \bn^*(\alpha_i)\ge s_ip_i,\quad i=1,...,n\}\ ,\nonumber\end{eqnarray}
\begin{eqnarray}{\mathcal M}^{br}_{0,n}(\Sigma,D,\bs,\bz,{\mathcal A})&=&\{[\bn:\PP^1\to
\Sigma,\bp]\in{\mathcal M}_{0,n}(\Sigma,D,\bs,\bz,{\mathcal A})\ :\ \nonumber\\
& &\quad\bn\ \text{is birational onto its image}\}\ ,\nonumber\end{eqnarray}
\begin{eqnarray}{\mathcal M}^{im}_{0,n}(\Sigma,D,\bs,\bz,{\mathcal A})&=&\{[\bn:\PP^1\to
\Sigma,\bp]\in{\mathcal M}^{br}_{0,n}(\Sigma,D,\bs,\bz,{\mathcal A})\ :\ \nonumber\\
& &\quad\bn\ \text{is an immersion}\}\ ,\nonumber
\end{eqnarray}
\begin{eqnarray}{\mathcal M}^{sing,1}_{0,n}(\Sigma,D,\bs,\bz,{\mathcal A})&=&\{[\bn:\PP^1\to
\Sigma,\bp]\in{\mathcal M}^{br}_{0,n}(\Sigma,D,\bs,\bz,{\mathcal A})\ :\ \nonumber\\
& &\quad\bn\ \text{is singular in}\ \PP^1\setminus\bp
\ \text{and smooth at}\ \bp\}\ ,\nonumber
\end{eqnarray}
\begin{eqnarray}{\mathcal M}^{sing,2}_{0,n}(\Sigma,D,\bs,\bz,{\mathcal A})&=&\{[\bn:\PP^1\to
\Sigma,\bp]\in{\mathcal M}^{br}_{0,n}(\Sigma,D,\bs,\bz,{\mathcal A})\ :\ \nonumber\\
& &\quad\bn\ \text{is singular at some point}\ p_i\in\bp\}\ .\nonumber
\end{eqnarray}

We shall consider the following strata in
$\Arc^{\sm}_s(\ring\Sigma^n)$:
\begin{enumerate}\item[(i)] The
subset $U^{im}(D)\subset \Arc^{\sm}_s(\ring\Sigma^n)$ is defined by the following conditions:

For any sequence $\bs=(s_1,...,s_n)\in\Z_{>0}^n$ summing up to $|\bs|\le s$ and for any
element $(\bz,{\mathcal A})\in U^{im}(D)$, where $\bz=(z_1,...,z_n)\in\ring\Sigma^n$,
${\mathcal A}=(\alpha_1,...,\alpha_n)$, $\alpha_i\in
\Arc_s(\Sigma,z_i)$, the family
${\mathcal M}_{0,n}(\Sigma,D,\bs,\bz,{\mathcal A})$
is empty if $|\bs|\ge-DK_\Sigma$, and is finite
if $|\bs|=-DK_\Sigma-1$.
Furthermore, in the latter case, 
all elements \mbox{$[\bn:\PP^1\to\Sigma,\bp]\in{\mathcal M}_{0,n}(\Sigma,D,\bs,\bz,{\mathcal A})$}
are represented by immersions $\bn:\PP^1\to\Sigma$ such that
$\bn^*(\alpha_i)=s_ip_i$, $1\le i\le n$. 
\item[(ii)] The subset $U^{im}_+(D)
\subset\Arc^{\sm}_s(\ring\Sigma^n)$ is defined by the following condition:

For any element $(\bz,{\mathcal A})\in U^{im}_+(D)$, there exists $\bs\in\Z_{>0}^n$ with $|\bs|
\ge-DK_\Sigma$ such that
${\mathcal M}^{im}_{0,n}(\Sigma,D,\bs,\bz,{\mathcal A})\ne\emptyset$.
\item[(iii)] The subset $U^{sing}_1(D)\subset\Arc^{\sm}_s(\ring\Sigma^n)$ is
defined by the following condition:

For any element $(\bz,{\mathcal A})\in U^{im}_+(D)$, there exists $\bs\in\Z_{>0}^n$ with $|\bs|=
-DK_\Sigma-1$ such that
${\mathcal M}^{sing,1}_{0,n}(\Sigma,D,\bs,\bz,{\mathcal A})\ne\emptyset$.
\item[(iv)] The subset $U^{sing}_2(D)\subset\Arc^{\sm}_s(\ring\Sigma^n)$ is
defined by the following condition:

For any element $(\bz,{\mathcal A})\in U^{sing}_2(D)$, there exists $\bs\in\Z_{>0}^n$ with
$|\bs|=-DK_\Sigma-1$ such that
${\mathcal M}^{sing,2}_{0,n}(\Sigma,D,\bs,\bz,{\mathcal A})\ne\emptyset$.
\item[(v)] The subset $U^{mt}(D)\subset\Arc^{\sm}_s(\ring\Sigma^n)$ is
defined by the following condition:

For any element $(\bz,{\mathcal A})\in U^{mt}(D)$, there exists $\bs\in\Z_{>0}^n$ with
$|\bs|=-DK_\Sigma-1$ and
$[\bn:\PP^1\to\Sigma,\bp]\in{\mathcal M}_{0,n}(\Sigma,D,\bs,\bz,{\mathcal A})$ such that $\bn$ is
a multiple cover of its image.
\end{enumerate}

\begin{proposition}\label{p4}
(1) The set $U^{im}(D)$ is Zariski open and dense in $\Arc^{\sm}_s(\ring\Sigma^n)$.

(2) If $U\subset U^{im}_+(D)$ is a component of codimension one in $\Arc^{\sm}_s(\ring\Sigma^n)$, then,
for a generic element $(\bz,{\mathcal A})\in U$ and any sequence $\bs\in\Z_{>0}^n$
with $|\bs|=-DK_\Sigma$, the set
${\mathcal M}^{im}_{0,n}(\Sigma,D,\bs,\bz,{\mathcal A})$ is either empty or finite,
and all of its elements $[\bn:\PP^1\to\Sigma,\bp]$ are presented by immersions and satisfy
$\bn^*(z_i)=s_ip_i$, $i=1,...,n$.

(3) If $U\subset U^{sing}_1(D)$ is a component of codimension one in
$\Arc^{\sm}_s(\ring\Sigma^n)$, then, for a generic element $(\bz,{\mathcal A})\in U$ and
any sequence $\bs\in\Z_{>0}^n$ with $|\bs|=-DK_\Sigma-1$, the set
${\mathcal M}^{sing,1}_{0,n}(\Sigma,D,\bs,\bz,{\mathcal A})$ is either empty or finite, whose
all elements $[\bn:\PP^1\to\Sigma,\bp]$ satisfy
$\bn^*(z_i)=s_ip_i$, $i=1,...,n$.

(4) If $U\subset U^{sing}_2(D)$ is a component of codimension one in
$\Arc^{\sm}_s(\ring\Sigma^n)$, then, for a generic element $(\bz,{\mathcal A})\in U$ and
any sequence $\bs\in\Z_{>0}^n$ with $|\bs|=-DK_\Sigma-1$, the set
${\mathcal M}^{sing,2}_{0,n}(\Sigma,D,\bs,\bz,{\mathcal A})$ is either empty or finite, whose
all elements $[\bn:\PP^1\to\Sigma,\bp]$ satisfy
$\bn^*(z_i)=s_ip_i$, $i=1,...,n$.

(5) If $U\subset U^{mt}(D)$ is a component of codimension one in
$\Arc^{\sm}_s(\ring\Sigma^n)$, then, for a generic element $(\bz,{\mathcal A})\in U$ and
any sequence $\bs\in\Z_{>0}^n$ with $|\bs|=-DK_\Sigma-1$, the following holds: Each element
$[\bn:\PP^1\to\Sigma,\bp]\in{\mathcal M}_{0,n}(\Sigma,D,\bs,\bz,{\mathcal A})$ satisfying
$C'=\bn(\PP^1)\in|D'|$, where $D=kD'$, $k\ge2$, admits a factorization
$$\bn:\PP^1\overset{\rho}{\longrightarrow}\PP^1\overset{\nu}{\longrightarrow} C'\hookrightarrow\Sigma$$
with $\rho$ a $k$-multiple ramified covering, $\nu$ the normalization, $\bp'=\rho(\bp)$, for
which one has
$$[\nu:\PP^1\to\Sigma,\bp']\in{\mathcal M}_{0,n}(\Sigma,D',\bs',\bz,{\mathcal A})\ ,$$
where $|\bs'|=-D'K_\Sigma$, and all branches $\nu\big|_{\PP^1,p'_i}$ are smooth.
\end{proposition}

{\bf Proof.} {\bf(1)} A general element of $[\bn:\PP^1\to\Sigma,\bp]\in{\mathcal M}_{0,n}(\Sigma,D)$
is represented by an immersion sending $\bp$ to $n$ distinct points of $\Sigma$
(cf. \cite[Lemma 9]{Itenberg_Kharlamov_Shustin:2012}).
Let $(\bz,{\mathcal A})\in\Arc^{\sm}_s(\ring\Sigma^n)$, and a sequence $\bs=(s_1,...,s_n)\in\Z_{>0}^n$
satisfy $|\bs|=-DK_\Sigma-1$. The fiber of the map
$\arc_{\bs}:{\mathcal M}_{0,n}(\Sigma,D)\to\prod_{i=1}^n\Arc^{\sm}_{s_i-1}(\Sigma)$, sending
an element $[\bn:\PP^1\to\Sigma,\bp]$ to the collection of arcs defined by the branches $\bn|_{\PP^1,
p_i}$, is either empty, or finite. Indeed, otherwise, by Lemma \ref{gus}(ii), we would get
a contradiction:
$$D^2\ge(D^2+DK_\Sigma+2)+|\bs|=D^2+1>D^2\ .$$ On the other hand,
$$\dim{\mathcal M}_{0,n}(\Sigma,D)=\dim\prod_{i=1}^n\Arc^{\sm}_{s_i-1}(\Sigma)=-DK_\Sigma-1+n\ ,$$
and hence the map $\arc_{\bs}$ is dominant. It follows, that, for a generic element
$(\bz,{\mathcal A})\in\Arc^{\sm}_s(\ring\Sigma^n)$ and any sequence $\bs\in\Z_{\ge0}^n$ such that
$|\bs|\le s$, one has: ${\mathcal M}^{im}_{0,n}(\Sigma,D,\bs,\bz,{\mathcal A})$ is empty if
$|\bs\ge-DK_\Sigma$, and ${\mathcal M}^{im}_{0,n}(\Sigma,D,\bs,\bz,{\mathcal A})$ is finite non-empty
if \mbox{$|\bs|=-DK_\Sigma-1$}. The same argument proves Claims (2) and (3) together with the fact that
$U^{im}_+(D)$ and $U^{sing}_1(D)$ have positive codimension in
$\Arc^{\sm}_s(\ring\Sigma^n)$.

Next, we will show that the sets $U^{sing}_2(D)$ and $U^{mt}(D)$ have positive codimension in
$\Arc^{\sm}_s(\ring\Sigma^n)$,
thereby completing the proof of
Claim (1), and will prove
Claims (4) and (5).

\smallskip

{\bf(2)} To proceed further, we introduce additional notations.
Let \mbox{$f:(\C,0)\to(C,z)\hookrightarrow(\Sigma,z)$} be the normalization of a reduced,
irreducible curve germ
$(C,z)$, and let $m_0,m_1,...$ be the multiplicities of $(C,z)$ and of its subsequent strict transforms
under blow-ups. We call this (infinite) sequence the {\it multiplicity sequence}
of $f:(\C,0)\to\Sigma$ and denote it $\overline m(f)$.
Note that, if $z_0=z$ and the infinitely near points
$z_1,...,z_j$, $0\le j\le s$, of $(C,z)$ belong to an arc from $\Arc^{\sm}_s(\Sigma,z)$,
then $m_0=...=m_{j-1}$ (see, for instance, \cite[Chapter III]{BK}). Such sequences $m_0,...,m_j$
will be called {\it smooth sequences}. Given smooth sequences
$\overline m_i=(m_{0i},...,m_{j(i),i})$ such
that $|\overline m_i|:=\sum_lm_{li}\le s$, $i=1,...,n$, denote by
${\mathcal M}_{0,n}(\Sigma,D,\{\overline m_i\}_{i=1}^n)$ the family of elements
$[\bn:\PP^1\to\Sigma,\bp]\in{\mathcal M}_{0,n}(\Sigma,D)$ such that $\bn$ is birational onto its image, and
$\overline m(\bn|_{\PP^1,p_i})$ contains $\overline m_i$ for every $i=1,...,n$.
Put
\begin{eqnarray}\Arc_s^{\sm}(\ring\Sigma^n,D,\{\overline m_i\}_{i=1}^n)&=&\{(\bz,{\mathcal A})
\in\Arc^{\sm}_s(\ring\Sigma^n)\ :\ \text{there exists}\nonumber\\
& &[\bn:\PP^1\to\Sigma,\bp]\in{\mathcal M}_{0,n}(\Sigma,D,\{\overline m_i\}_{i=1}^n)\nonumber\\
& &\text{such that}\ \bn(\bp)=\bz\ \text{and}\ \bn^*(\alpha_i)\ge|\overline m_i|p_i,\ i=1,...,n\}
\nonumber\end{eqnarray}

\smallskip

{\bf(3)} We now prove Claim (4) together with the fact that $U^{sing}_2(D)$ has positive codimension
in $\Arc^{\sm}_s(\ring\Sigma^n)$.

Let $(\bz,{\mathcal A})$ be a generic element of a top-dimensional
component $U\subset U^{sing}_2(D)$, a sequence $\bs\in\Z_{>0}^n$ satisfy $|\bs|=-DK_\Sigma-1$, and
$[\bn:\PP^1\to\Sigma,\bp]\in{\mathcal M}^{br}_{0,n}(\Sigma,D,\bs,\bz,{\mathcal A})$
have singular branches among
$\bn|_{\PP^1,p_i}$, $i=1,...,n$. Let $\overline m_i=(m_{0i},...,m_{j(i),0})$ be a
smooth multiplicity sequence of
the branch $\bn|_{\PP^1,p_i}$ such that $|\overline m_i|\ge s_i$. Denote by
${\mathcal V}$ the germ at $[\bn:\PP^1\to\Sigma,\bp]$ of a top-dimensional component
of ${\mathcal M}_{0,n}(\Sigma,D,\{\overline m_i\}_{i=1}^n)$. Without loss of generality,
we can suppose that
${\mathcal M}^{br}_{0,n}(\Sigma,D,\bs,\bz,{\mathcal A})\subset{\mathcal M}_{0,n}(\Sigma,D,
\{\overline m_i\}_{i=1}^n)$
and $U\subset\Arc_s^{\sm}(\ring\Sigma^n,D,\{\overline m_i\}_{i=1}^n)$.

Note that $[\bn:\PP^1\to\Sigma,\bp]$ is isolated in
${\mathcal M}^{br}_{0,n}(\Sigma,D,\bs,\bz,{\mathcal A})$. Indeed, otherwise Lemma \ref{gus}(ii) would
yield a contradiction:
$$D^2\ge(D^2+DK_\Sigma+2)+\sum_{i=1}^n(m_{0i}-1+|\overline m_i|)\ge(D^2+DK_\Sigma+2)+|\bs|=
D^2+1>D^2\ .$$

Next, we can suppose that $m_{0i}\ge2$ as $1\le i\le r$ for some $1\le r\le n$, and that
$m_{0i}=1$ for $r<i\le n$.

Consider the case when $|\overline m_i|=s_i$ for all $i=1,...,n$. We claim that then
\begin{equation}
\dim {\mathcal V}\le\sum_{i=1}^nj(i)+n+r-1\ .\label{ee1}\end{equation}
If so, we would get
$$\dim U\le\sum_{i=1}^n(s-j(i))+n-r+\dim {\mathcal V}
\le n(s+2)-1=\dim\Arc^{\sm}_s(\ring\Sigma^n)-1\ ,$$
and the equality would yield $(\bn')^*(\bz,{\mathcal A})=\sum_{i=1}^ns_i=-DK_\Sigma-1$ for each element
\mbox{$[\bn':\PP^1\to\Sigma,\bp']\in{\mathcal M}^{sing,2}_{0,n}(\Sigma,D,\bs,\bz,{\mathcal A})$}
with generic $(\bz,{\mathcal A})\in U$, as required in Claim (3).
To prove (\ref{ee1}), we show that the
assumption
\begin{equation}\dim{\mathcal V}\ge\sum_{i=1}^nj(i)+n+r\label{ee2}
\end{equation}
leads to contradiction. Namely,
we impose $\sum_{i=1}^nj(i)+n+r-1$
conditions, defining a positive-dimensional subfamily of ${\mathcal V}$
containing $[\bn:\PP^1\to\Sigma,\bp]$,
and apply Lemma \ref{gus}. It is enough to consider the following situations:
\begin{enumerate}\item[(a)] $1\le r<n$;
\item[(b)] $1<r=n$, $j(1)>0$;
\item[(c)] $1=r=n$, $j(1)>0$, $m_{01}>m_{j(1),1}$;
\item[(d)] $r=n$, $j(1)=...=j(n)=0$;
\item[(e)] $1=r=n$, $j(1)>0$, $m_{01}=...=m_{j(1),1}$.
\end{enumerate}
In case (a), we
fix the position of $z_i$ and of the next $j(i)$ infinitely near points for $i=1,...,r$, and
the position of additional
$\sum_{i=r+1}^nj(i)+n-r-1$ smooth
points on $C=\bn(\PP^1)$, obtaining a positive-dimensional subfamily of $U$
and a contradiction by Lemma \ref{gus}:
\begin{eqnarray}D^2&\ge&(D^2+DK_\Sigma+2)+\sum_{i=1}^r(m_{0i}-1+|\overline m_i|)+
\sum_{i=r+1}^nj(i)+n-r-1\nonumber\\
&=&D^2+\sum_{i=1}^r(m_{0i}-1)>D^2\ .\nonumber\end{eqnarray}
In case (b), we fix the position of $\bz$ and of additional infinitely near points:
$j(1)-1$ points for $z_1$, and $j(i)$ points for all $2\le i\le n$. These conditions define
a positive-dimensional
subfamily of $U$, which implies a contradiction by Lemma \ref{gus}:
\begin{eqnarray}D^2&\ge&(D^2+DK_\Sigma+2)+\sum_{i=2}^r(m_{0i}-1+|\overline m_i|
)+(m_{01}-1)+|\overline m_1|-m_{j(1),1}
\nonumber\\
&\ge&D^2+\sum_{i=2}^n(m_{0i}-1)>D^2\ .\nonumber\end{eqnarray}
In case (c), the same construction similarly leads to a contradiction:
$$D^2\ge(D^2+DK_\Sigma+2)+(m_{01}-1)+\sum_{0\le k<j(1)}m_{k1}\ge
(D^2+DK_\Sigma+2)+|\overline m_1|=D^2+1>D^2\ .$$
In case (d), we fix the position of $z_i$, $1<i\le n$, and of one more smooth point of $C=\bn(\PP^1)$.
Thus, Lemma \ref{gus}, applied to the obtained positive-dimensional family, yields a contradiction:
$$D^2\ge(D^2+DK_\Sigma+2)+\sum_{i=1}^n(m_{0i}-1)+\sum_{1<i\le n}m_{0i}+1
=D^2+\sum_{1<i\le n}(m_{0i}-1)+1>D^2\ .$$
In case (e), relation (\ref{ee2}) reads $\dim{\mathcal V}\ge j(1)+2=\dim\Arc_{j(1)}(\Sigma)$.
As noticed above,
the map $\arc_{j(1)}:{\mathcal V}\to\Arc_{j(1)}(\Sigma)$ is finite. Hence, $\dim{\mathcal V}=j(1)+2$, and
(due to the general choice of $\xi=[\bn:\PP^1\to\Sigma,p]\in{\mathcal V}$) the germ $({\mathcal V},\xi)$
diffeomorphically maps onto the germ of $\Arc_{j(1)}(\Sigma)$ at $\pi(\xi)$. Observe that the fragment
$(m_{01},...,m_{j(1),1},m_{j(1)+1,1})$ of the multiplicity sequence of $\bn|_{\PP^1,p}$ is a smooth
sequence.
That means, the map of $({\mathcal V},\xi)$ to
$\Arc_{j(1)+1}(\Sigma)$ defines a section
$\sigma:(\Arc_{j(1)}(\Sigma),\pi(\xi))\to\Arc_{j(1)+1}(\Sigma)$, satisfying the hypotheses of Lemma
\ref{arc}. So, we take the curve $\Lambda$, defined in Lemma \ref{arc}, and apply
Lemma \ref{gus}(iii):
\begin{eqnarray}D^2&\ge&(D^2+DK_\Sigma+2)+(m_{01}+...+m_{j(1),1}+m_{j(1)+1,1})-1\nonumber\\
&\ge&(D^2+DK_\Sigma+2)+|\overline m_1|=D^2+1>D^2\ ,\nonumber\end{eqnarray} which completes the proof of
(\ref{ee1}).

Consider the case when $\sum_{i=1}^n|\overline m_i|>-DK_\Sigma-1$ and show that then
$\dim U\le\dim\Arc^{\sm}_s(\ring\Sigma^n)-2$. The preceding consideration reduces the problem
to the case
$$r=n\quad\text{and}
\quad\sum_{i=1}^n|\overline m_i|-m_{j(n),n}<-DK_\Sigma-1<\sum_{i=1}^n|\overline m_i|\ ,$$ in which we
need to prove that
\begin{equation}\dim{\mathcal V}\le\sum_{i=1}^nj(i)+2n-2\ .\label{ee3}\end{equation}
We assume that
\begin{equation}\dim{\mathcal V}\ge\sum_{i=1}^nj(i)+2n-1\label{ee4}\end{equation}
and derive a contradiction in the same manner as for
(\ref{ee2}). We shall separately treat several possibilities:
\begin{enumerate}\item[(a)] $j(n)=0$;
\item[(b)] $j(n)>0$.
\end{enumerate}
In case (a), we fix the position of $z_i$ and of the additional $j(i)$ infinitely near points
for all $i=1,...,n-1$, thereby cutting off ${\mathcal V}$
a positive-dimensional subfamily, and hence by Lemma
\ref{gus} we get a contradiction:
\begin{eqnarray}D^2&\ge&(D^2+DK_\Sigma+2)+\sum_{i=1}^{n-1}(m_{0i}-1+|\overline m_i|)+m_{0n}-1\nonumber\\
&\ge&(D^2+DK_\Sigma+2)+\sum_{i=1}^n|\overline m_i|-1\ge D^2+1>D^2\ .\nonumber\end{eqnarray}
In case (b), we again fix the position of $z_i$ and of the additional $j(i)$ infinitely near points
for all $i=1,...,n-1$, thereby cutting off ${\mathcal V}$
a subfamily ${\mathcal V}'$ of dimension $\ge j(n)+1$. Consider the map
$\arc_{j(n)-1}:{\mathcal V}'\to\Arc_{j(n)-1}(\Sigma)$ and note that
$\dim\Arc_{j(n)-1}(\Sigma)=j(n)+1\le\dim{\mathcal V}'$. If $\dim\pi(
{\mathcal V}')\le j(n)$, fixing the position of
$z_n$ and of $j(n)-1$ additional infinitely near points, we
obtain a positive-dimensional subfamily of ${\mathcal V}'$, and hence a
contradiction by Lemma \ref{gus}:
\begin{eqnarray}D^2&\ge&(D^2+DK_\Sigma+2)+\sum_{i=1}^{n-1}(m_{0i}-1+|\overline m_i|)+
(m_{0n}-1)+|\overline m_n|-m_{j(n),n}\nonumber\\
&\ge&(D^2+DK_\Sigma+2)+\sum_{i=1}^n|\overline m_i|-1\ge D^2+1>D^2\ .\nonumber\end{eqnarray}
If $\dim\pi({\mathcal V}')=j(n)+1$, the preceding argument yields that \mbox{$\dim
{\mathcal V}'=j(n)+1$}, and we can
suppose that the germ of ${\mathcal V}'$ at the initially chosen element \mbox{$\xi=[\bn:\PP^1\to
\Sigma,\bp]\in{\mathcal V}$} is diffeomorphically mapped onto the germ of
$\Arc_{j(n)-1}(\Sigma)$ at $\arc_{j(n)-1}(\xi)$. Thus, we obtain a section
$\sigma:(\Arc_{j(n)-1}(\Sigma),\pi(\xi))\to\Arc_{j(n)}(\Sigma)$ defined by the map
$({\mathcal V}',\xi)\to\Arc_{j(n)}(\Sigma)$. It satisfies the hypotheses of Lemma \ref{arc},
which allows one to construct a smooth curve $\Lambda$ as in Lemma \ref{arc}
and apply Lemma \ref{gus}(iii):
$$D^2\ge(D^2+DK_\Sigma+2)+\sum_{i=1}^{n-1}(m_{0i}-1+|\overline m_i|)+|\overline m_n|-1
\ge D^2+1>D^2\ ,$$ a contradiction.

The proof of Claim (4) is completed.

\smallskip

{\bf(4)} It remains to consider the set $U^{mt}(D)$.
Let $(\bz,{\mathcal A})\in U^{mt}(D)$, $\bs\in\Z_{>0}^n$ satisfy $|\bs|=-DK_\Sigma-1$, and
$[\bn:\PP^1\to\Sigma,\bp]\in{\mathcal M}_{0,n}(\Sigma,D,\bs,\bz,{\mathcal A})$ be such that
$\bn$ is a $k$-multiple (ramified) covering of its image $C=\bn(\PP^1)$, $k\ge2$.
We have $C\in|D'|$, where
$kD'=D$, and $\nu^*(\alpha_i)\ge s'_ip'_i$, $\rho^*(p'_i)\ge l_ip_i$,
where $l_is'_i\ge s_i$ for all $i=1,...,n$. Since $l_i\le k$ for all $i=1,...,n$, it follows that
$$\sum_{i=1}^ns'_i\ge\frac{|\bs|}{k}=\frac{-DK_\Sigma-1}{k}=-D'K_\Sigma-\frac{1}{k}>-D'K_\Sigma-1\ .$$
This yields that
$U^{mt}(D)$ has positive codimension in $\Arc^{\sm}_s(\ring\Sigma^n)$, and, furthermore, if not all
branches $\nu\big|_{\PP^1,p'_i}$, $i=1,...,n$, are smooth, the codimension of $U^{mt}(D)$
in $\Arc^{\sm}_s(\ring\Sigma^n)$ is at least $2$.
The proof of Claim (4) and thereby of Claim (1) is completed.
\proofend

\subsection{Families of curves and arcs on generic del Pezzo surfaces}\label{secGDP}
Let $\Sigma$ be a smooth del Pezzo surface of degree $1$ satisfying the following condition:

(GDP) There are only finitely many effective divisor classes $D\in\Pic(\Sigma)$ satisfying $-DK_\Sigma=1$,
and for any such divisor $D$, the linear system $|D|$ contains only finitely many rational curves,
all these rational curves are immersed, and any two curves $C_1\ne C_2$ among them intersect
in $C_1C_2$ distinct points.

By \cite[Lemmas 9 and 10]{Itenberg_Kharlamov_Shustin:2012}, these del Pezzo surfaces
form an open dense subset in the space of del Pezzo surfaces of degree $1$.

Let us fix an effective divisor $D\in\Pic(\Sigma)$ such that $-DK_\Sigma-1\ge3$.

\begin{proposition}\label{p5}
In the notation of Section \ref{sec2.2}, let $(\bz_0,{\mathcal A}_0)$ be a generic element of a component
$U$ of $U^{mt}(D)$ having codimension one in $\Arc^{\sm}_s(\ring\Sigma^n)$,
a sequence $\bs\in\Z_{>0}^n$ satisfy $|\bs|=-DK_\Sigma-1$, and $[\bn_0:\PP^1\to\Sigma,\bp_0]\in
{\mathcal M}_{0,n}(\Sigma,D,\bs,\bz_0,{\mathcal A}_0)$ be such that $\bn_0$ covers its image with
multiplicity
$k\ge2$ so that $\bn_0(\PP^1)\in|D'|$, where $D=kD'$, and
$\bn_0=\nu\circ\rho$ with $\nu:\PP^1\to C'$ the normalization,
$\rho:\PP^1\to\PP^1$ a $k$-fold ramified covering. Assume that $(\bz_t,{\mathcal A}_t)\in
\Arc^{\sm}_s(\ring\Sigma^n)$, $t\in(\C,0)$, is the germ at $(\bz_0,{\mathcal A}_0)$ of a
generic one-dimensional family such that $(\bz_t,{\mathcal A}_t)\not\in U^{mt}(D)$ as $t\ne0$,
and assume that
there exists a family $[\bn_t:\PP^1\to\Sigma,\bp_t]\in{\mathcal M}_{0,n}(\Sigma,D,\bs,\bz_t,{\mathcal
A}_t)$ extending the element $[\bn_0:\PP^1\to\Sigma,\bp_0]$. Then 
$n=3$, $k=2$, $-D'K_\Sigma=3$, $\bs=(2,2,1)$, and
$[\nu:\PP^1\to C'\hookrightarrow\Sigma,\bp'_0]\in{\mathcal M}_{0,3}(\Sigma,D',\bs',\bz_0,{\mathcal A}_0)$,
where $\bp'=\rho(\bp_0)$ and $\bs'=(1,1,1)$. Furthermore, the family
$[\bn_t:\PP^1\to\Sigma,\bp_t]$, $t\in(\C,0)$, is smooth and isomorphically projects onto the
family $(\bz_t,{\mathcal A}_t)$, $t\in(\C,0)$.
\end{proposition}

{\bf Proof.}
Note, first, that by the assumption (GDP) and Proposition \ref{p4}(2,5), the map $\bn_0:\PP^1\to\Sigma$
is an immersion,
and (in the notation of Proposition \ref{p4}(5))
\begin{equation}\nu^*(\alpha_i)=s'_ip'_i,\ i=1,...,n,\quad \sum_{i=1}^ns'_i=-D'K_\Sigma\ .\label{emt2}
\end{equation}
Furthermore, if $C'=\bn_0(\PP^1)\in|D'|$, where $D=kD'$, then $(D')^2>0$, since the assumption
$-DK_\Sigma\ge4$ yields $D^2\ge2$ by the adjunction formula. Hence, in the deformation $\bn_t:\PP^1\to
\Sigma$, $t\in(\C.0)$,
in a neighborhood of each singular point $z$ of $C'$, there appear singular points of $C_t=\bn_t(\PP^1)$,
$t\ne0$,
with total $\delta$-invariant at least $k^2\delta(C',z)$, which implies
\begin{equation}k^2\left(\frac{(D')^2+D'K_\Sigma}{2}+1\right)\le\frac{k^2(D')^2+kD'K_\Sigma}{2}+1\ ,
\label{emt0}\end{equation} and hence
\begin{equation}-D'K_\Sigma\ge\frac{2k+2}{k}\quad\text{or, equivalently,}\quad-D'K_\Sigma\ge3\ .
\label{emt1}\end{equation}

Let $\rho^*(p'_i)\ge l_ip_i$, $i=1,...,n$. We can suppose that $k\ge l_1\ge...\ge l_n$. Then
\begin{equation}\sum_{i=1}l_is'_i\ge-kD'K_\Sigma-1\quad\Longrightarrow\quad
\sum_{i=1}^n(l_i-1)s'_i\ge-(k-1)D'K_\Sigma-1\ .
\label{emt3}\end{equation}
If $l_1\le k-1$, then (\ref{emt2}) and (\ref{emt3}) yield
$$-(k-2)D'\ge-(k-1)D'K_\Sigma-1\quad\Longrightarrow\quad-D'K_\Sigma\le 1\ ,$$ forbidden by (\ref{emt1}),
and hence
\begin{equation}l_1=k\ .\label{emt4}\end{equation} By Riemann-Hurwitz, $\sum_{i>1}(l_i-1)\le k-1$,
and then it follows
from (\ref{emt3}) that
\begin{equation}(k-1)(-D'K_\Sigma-(n-1))+(k-1)\ge-(k-1)D'K_\Sigma-1\ ,\label{emt5}
\end{equation} or, equivalently
\begin{equation}(n-2)(k-1)\le1\ ,\label{emt6}\end{equation}
which in view of Riemann-Hurwitz and (\ref{emt4})-(\ref{emt6}) leaves the following options:
\begin{itemize}\item either $n=1$,
\item or $n=2$, $\bs=(k(-D'K_\Sigma-1),(k-1))$,
\item or $n=2$, $\bs=(ks'_1,ks'_2)$, $s'_1+s'_2=-D'K_\Sigma$,
\item or $n=3$, $\bs=(2(-D'K_\Sigma-2),2,1)$.
\end{itemize}
Let us show that $s'_1>1$ is not possible. Indeed, otherwise, in suitable local coordinates $x,y$ in
a neighborhood of $z_1$ in $\Sigma$,, we would have $z_1=(0,0)$, $C'=\{y=0\}$, $\bn_0:(\PP^1,p_1)\to
(\Sigma,z_1)$ acts by $\tau\in(\C,0)\simeq(\PP^1,p_1)\mapsto(\tau^k,\tau)$, and we also may assume that
the family of arcs $\alpha_{1,t}$ is centered at $z_1$ and given by
$y=\sum_{i\ge s'_1}a_i(t)x^i$ with $a_i(0)\ne0$, $i\ge s'_1$. Then $\bn_t:(\PP^1,p_{1,t})\to(\Sigma,
z_1)$ can be expressed via $\tau\in(\C,0)\simeq(\PP^1,p_{1,t})\mapsto(\tau^k+tf(t,\tau),tg(t,\tau))$, which
contradicts the requirement $\bn_t^*(\alpha_{1,t})\ge(ks'_1-1)p_{1,t}$ equivalently written as
$$t\cdot g(t,\tau)\equiv\sum_{i\ge s'_1}a_i(t)(\tau^k+tf(t,\tau))^i\mod (\tau^k+tf(t,\tau))^{ks'_1-1}\ ,$$ since
the term $a_{s'_1}(0)\tau^{ks'_1}$ does not cancel out here in view of $k\ge2$.

Thus, in view of (\ref{emt1}), we are left with $n=3$, $k=2$, $\bs'=(1,1,1)$, and $\bs=(2,2,1)$.
Without loss of generality, for $(z_t,{\mathcal A}_t)$, $t\in(\C,0)$, we can choose the family consisting
of two fixed points $z_{1,0},z_{2,0}$ and fixed arcs $\alpha_{1,0},\alpha_{2,0}$
(transversal to $C'$), and of
a point $z_{3,\tau}$ mowing along the germ $\Lambda$ of a smooth curve transversally intersecting $C'$
at $z_{3,0}$ ($\tau$ being a regular parameter on $\Lambda$). We then claim that the evaluation
$$[\bn_t:\PP^1\to\Sigma,\bp_t]\mapsto\bn_t(p_{3.t})=z_{3,\tau(t)}$$ is one-to-one, completing the
proof of Proposition \ref{p5}. So, we establish the formulated claim arguing on the contrary: if
some point $z_{3,\tau}$, $\tau\ne0$, has two preimages, then the curves $C_1=\bn_{t_1}(\PP^1)$,
$C_2=\bn_{t_2}(\PP^1)$ intersect with total multiplicity $\ge5$ at $z_{1,0},z_{2.0},z_{3,\tau}$,
and intersect with multiplicity $\ge\delta(C',z)$ in a neighborhood of each point $z\in\Sing(C')$,
which altogether leads to a contradiction:
$$C_1C_2\ge5+4((D')^2+D'K_\Sigma+2)=5+D^2-4=D^2+1\ .\quad\text{\proofend}$$

The compactification $\overline{\mathcal M}_{0,n}(\Sigma,
D,\bs,\bz,{\mathcal A})$ of the space ${\mathcal M}_{0,n}(\Sigma,
D,\bs,\bz,{\mathcal A})$ is obtained by adding the elements
$[\bn:\widehat C\to\Sigma,\bp]$, where \begin{itemize}\item
$\widehat C$ is a tree formed by $k\ge2$ components
$\widehat C^{(1)},...,\widehat C^{(k)}$ isomorphic to $\PP^1$;
\item the points of $\bp$ are distinct but allowed to be at the nodes of $\widetilde C$;
\item $[\bn:\widehat C^{(j)}\to\Sigma,\widehat C^{(j)}\cap\bp]\in{\mathcal M}_{0,
|C^{(j)}\cap\bp|}(\Sigma,D^{(j)},\bs^{(j)},
\bz,{\mathcal A})$, where we suppose that the integer vector $\bs^{(j)}\in\Z_{\ge0}^n$ has coordinates
$s^{(j)}_i>0$ or $s^{(j)}_i=0$ according as $p_i$ belongs to $\widehat C^{(j)}$ or not,
$j=1,...,k$;
\item $\sum_{j=1}^kD^{(j)}=D$, where $D^{(j)}\ne0$, $j=1,...,k$, and $\sum_{j=1}^k\bs^{(j)}=\bs$.
\end{itemize}
One can view this compactification as the image of the closure of 
${\mathcal M}_{0,n}(\Sigma,
D,\bs,\bz,{\mathcal A})$ in the moduli space 
of stable maps $\overline{\mathcal M}_{0,n}(\Sigma,D)$ under the morphism,
which contracts the components of the source curve that are mapped to points.
Notice that in our compactification the source curves $\widehat C$ may be not nodal, and the
marked points may appear at intersection points of components of a (reducible) source curve.

Introduce the set
$U^{red}(D)\subset\Arc^{\sm}_s(\ring\Sigma^n)$ defined by the following condition:
For any element $(\bz,{\mathcal A})\in U^{red}(D)$, there exists $\bs\in\Z_{>0}^n$ with
$|\bs|\ge-DK_\Sigma-1$ such that
$\overline{\mathcal M}_{0,n}(\Sigma,D,\bs,\bz,{\mathcal A})\setminus{\mathcal M}_{0,n}
(\Sigma,D,\bs,\bz,{\mathcal A})\ne\emptyset$.

\begin{proposition}\label{p6}
The set $U^{red}(D)$ has positive codimension in $\Arc^{\sm}_s(\ring\Sigma^n)$.
Let $(\bz,{\mathcal A})$ be a generic element of
a component of $U^{red}(D)$ having codimension one in $\Arc^{\sm}_s(\ring\Sigma^n)$, and let
$(\bz_t,{\mathcal A}_t)\in\Arc^{\sm}_s(\ring\Sigma^n)$, $t\in(\C,0)$, be a generic family
which transversally intersects $U^{red}(D)$ at $(\bz_0,{\mathcal A}_0)=(\bz,{\mathcal A})$.

(1) Given any vector $\bs\in\Z_{>0}^n$ such that $|\bs|=-DK_\Sigma-1$, the set
$\overline{\mathcal M}_{0,n}(\Sigma,D,\bs,\bz,{\mathcal A})\setminus{\mathcal M}_{0,n}
(\Sigma,D,\bs,\bz,{\mathcal A})$ is
either empty, or finite. Moreover, let
$$[\bn:\widehat C\to\Sigma,\bp]\in\overline{\mathcal M}_{0,n}(\Sigma,D,\bs,\bz,{\mathcal A})
\setminus{\mathcal M}_{0,n}(\Sigma,D,\bs,\bz,{\mathcal A})$$ extend to a family
\begin{equation}[\bn_\tau:\widehat C_\tau\to\Sigma,\bp_\tau]\in\overline{\mathcal M}_{0,n}
(\Sigma,D,\bs,\bz_{\varphi(\tau)},{\mathcal A}_{\varphi(\tau)}),\quad\tau\in(\C,0)\ ,
\label{ered6}\end{equation} for some morphism
$\varphi:(\C,0)\to(\C,0)$. Then $[\bn:\widehat C\to\Sigma,\bp]$ is as follows:
\begin{enumerate}\item[(1i)] either $\widehat C=C^{(1)}\cup C^{(2)}$, where $C^{(1)}\simeq C^{(2)}
\simeq\PP^1$,
$\bn(C^{(1)})\ne\bn(C^{(2)})$, and
\begin{itemize}\item
the map $\bn:\widehat C^{(j)}\to\Sigma$ is an immersion and $\bz\cap \Sing(C^{(j)})=\emptyset$ for $j=1,2$,
\item $|\bp\cap\widehat C^{(1)}\cap\widehat C^{(2)}|\le1$,
\item $[\bn:\widehat C^{(j)}\to\Sigma,\bp\cap\widehat C^{(j)}]\in
{\mathcal M}_{0,|\bp\cap\widehat C^{(j)}|}(\Sigma,D^{(1)},\bs^{(j)},\bz,{\mathcal A})$, $j=1,2$,
where $D^{(1)}+D^{(2)}=D$, $\bs^{(1)}+\bs^{(2)}=\bs$, $|\bs^{(1)}|=-D^{(1)}K_\Sigma$,
$|\bs^{(2)}|=-D^{(2)}K_\Sigma-1$, and, moreover,
$(\bn|_{C^{(j)}})^*({\mathcal A})=\sum_{i=1}^ns^{(j)}_ip_i$ for $j=1,2$;
\end{itemize}
\item[(1ii)] or $n=1$, $\bz=z_1\in\Sigma$, ${\mathcal A}=\alpha_1\in\Arc^\sm_s(\Sigma,z)$,
$\bp=p_1\in\widehat C$,
$D=kD'$, where $k\ge2$ and $-D'K_\Sigma\ge3$, and the following holds
\begin{itemize} \item $\widehat C$ consists of few components having $p_1$ as a common point, and each of them is mapped onto the same
immersed rational curve $C\in|D'|$;
\item $z_1$ is a smooth point of $C$, and $(C\cdot\alpha_1)=-D'K_\Sigma$.
\end{itemize}
\item[(1iii)] or $D=kD'+D''$, where $k\ge2$, $-D'K_\Sigma\ge2$, $D''\ne0$, $\widetilde C=
\widetilde C'\cup...\cup\widetilde C''$, where
\begin{itemize}\item $\widehat C'\simeq\PP^1$, $\bn:\widehat C''\to CX''\hookrightarrow\Sigma$ is an immersion, where $C''\in|D''|$,
\item the components of $\widehat C'$ have a common point $p_1$ and are disjoint from
$p_2,...,p_n$, and each of them is mapped onto the same
immersed rational curve $C'\in|D'|$,
\item $z_1$ is a smooth point of $C'$, and $(C'\cdot\alpha_1)=-D'K_\Sigma$.
\end{itemize}
\end{enumerate}

(2) In case (1i), \begin{itemize}\item if $\bp\cap\widehat C^{(1)}\cap\widehat C^{(2)}=\emptyset$, there is
a unique family of type (\ref{ered6}), and it is smooth, parameterized by $\tau=t$;
\item if $\widehat C^{(1)}\cap\widehat C^{(2)}=\{p_1\}$, then
there are precisely $\kappa=\min\{s^{(1)}_1,s^{(2)}_1\}$ families of type (\ref{ered6}), and
for each of them $t=
\tau^{\kappa/d}$, where $d=\text{gcd}(s^{(1)}_1,s^{(2)}_1)$.\end{itemize}
\end{proposition}

{\bf Proof.} If $[\bn:\widehat C\to\Sigma,\bp]\in\overline{\mathcal M}_{0,n}(\Sigma,D,\bs,\bz,
{\mathcal A})$ with a
generic $(\bz,{\mathcal A})\in\Arc^{\sm}_s(\ring\Sigma^n)$ and $\widehat C$ consisting of $m\ge1$
components, then by Propositions
\ref{p4} and \ref{p5} one obtains $m=1$ and $\bn$ immersion. Hence, $U^{red}(D)$ has positive
codimension in $\Arc^{\sm}_s(\ring\Sigma^n)$.
Suppose that $(\bz,{\mathcal A})$ satisfies the hypotheses of Proposition. Then the finiteness of
$\overline{\mathcal M}_{0,n}(\Sigma,D,\bs,\bz,{\mathcal A})\setminus{\mathcal M}_{0,n}(\Sigma,D,
\bs,\bz,{\mathcal A})$ and the asserted structure of its elements
follows from Propositions \ref{p4} and \ref{p5}, provided we show that
\begin{enumerate}
\item[(a)] there are no two components $\widehat C',\widehat C''$ of $\widehat C$ such that
$\bn(\widehat C')\ne\bn(\widehat C'')$, $\bn_*(\widehat C')\in|D'|$, $\bn_*(\widehat C'')\in|D''|$, and
$\deg(\bn|_{\widehat C'})^*{\mathcal A}\ge-D'K_\Sigma$, $\deg(\bn|_{\widehat C''})^*{\mathcal A}
\ge-D''K_\Sigma$,
\item[(b)] in cases (1ii) and (1iii) we have inequalities $-D'K_\Sigma\ge3$ and
$-D'K_\Sigma\ge2$, respectively.
\end{enumerate}
The proof of Claim (a) can easily be reduced to the case when $\bn|_{\widehat C'}$ and $\bn|_{\widehat
C''}$ are immersions, and $\deg(\bn|_{\widehat C'})^*\alpha_1=-D'K_\Sigma=\deg(\bn|_{\widehat C''})^*
\alpha_1=-D''K_\Sigma$. However, in such a case, the dimension and generality assumptions
yield that there exists the germ at $C''$ of the family of rational curves $C''_t\in|D''|$,
$t\in(\C,0)$,
such that $(C''_t\cdot C')_{y_t}\ge-D''K_\Sigma$ for some family of points $y_t\in(C',z_1)$,
$t\in(\C,0)$, which together with Lemma \ref{gus}(iii) implies a contradiction:
$$(D'')^2\ge((D'')^2+D''K_\Sigma+2)+(-D''K_\Sigma-1)=(D'')^2+1\ .$$
Claim (b) in the case (1ii) follows from inequalities (\ref{emt0}) and (\ref{emt1}). In case (1iii), we
perform similar estimations. If the curves $C'$ and $C''$ intersect at $z_1$, then $(C'\cdot C'')_{z_1}=
\min\{-D'K_\Sigma,-D''K_\Sigma-1\}$, and we obtain
$$\frac{(kD'+D'')^2+(kD'+D'')K_\Sigma}{2}+1\ge k^2\left(\frac{(D')^2+D'K_\Sigma}{2}+1\right)$$
$$+k(D'D''-(C'\cdot C'')_{z_1})+\frac{(D'')^2+D''K_\Sigma}{2}+1$$
$$\Longleftrightarrow\quad\begin{cases}(k-1)(-D'K_\Sigma)+2(-D''K_\Sigma-1)\ge2k,\quad&\text{if}\ -D'K_\Sigma\ge-D''K_\Sigma-1,\\
(k+1)(-D'K_\Sigma)\ge2k,\quad & \text{if}\ -D'K_\Sigma\le-D''K_\Sigma-1\end{cases}$$
$$\Longrightarrow\quad-D'K_\Sigma\ge2\ .$$
If the curves $C'$ and $C''$ do not meet at $z_1$, then we obtain
$$\frac{(kD'+D'')^2+(kD'+D'')K_\Sigma}{2}+1\ge k^2\left(\frac{(D')^2+D'K_\Sigma}{2}+1\right)$$
$$+k(D'D''-1)+\frac{(D'')^2+D''K_\Sigma}{2}+1\quad\Longleftrightarrow\quad-D'K_\Sigma\ge2\ .$$


\smallskip

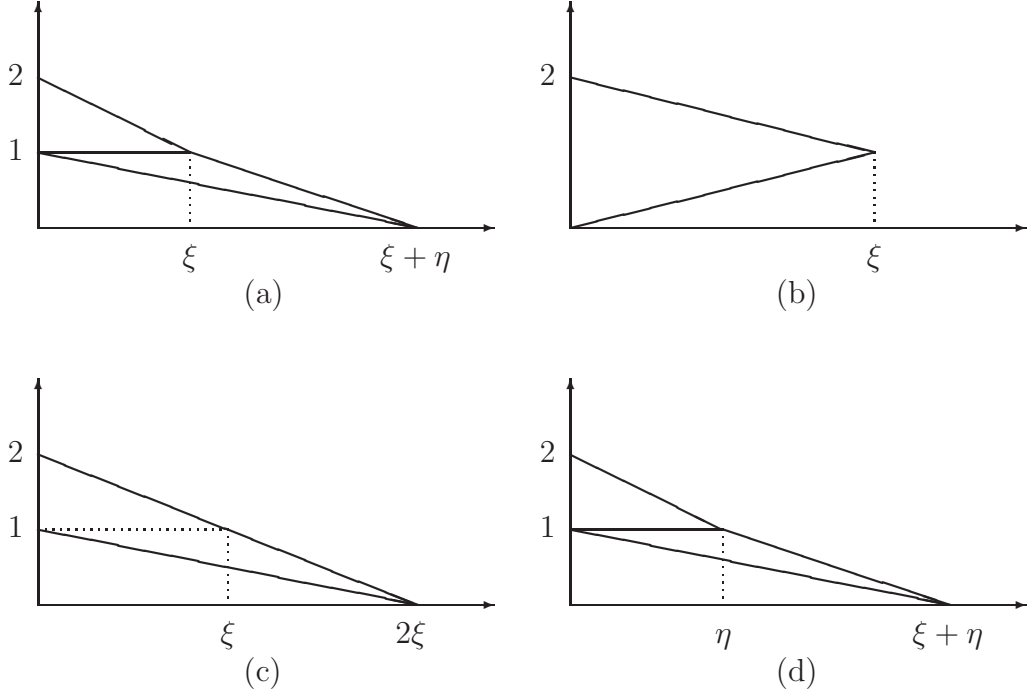
\begin{figure}
\setlength{\unitlength}{1cm}
\begin{picture}(14,9)(0,0)
\thinlines\put(1,1){\vector(1,0){6}}\put(8,1){\vector(1,0){6}}
\put(1,1){\vector(0,1){3}}\put(8,1){\vector(0,1){3}}
\put(1,6){\vector(1,0){6}}\put(8,6){\vector(1,0){6}}
\put(1,6){\vector(0,1){3}}\put(8,6){\vector(0,1){3}}
\dottedline{0.1}(1,2)(3.5,2)\dottedline{0.1}(3.5,2)(3.5,1)
\dottedline{0.1}(10,2)(10,1)\dottedline{0.1}(3,7)(3,6)
\dottedline{0.1}(12,7)(12,6)
\thicklines \put(1,2){\line(5,-1){5}}\put(1,3){\line(5,-2){5}}
\put(8,2){\line(5,-1){5}}\put(8,2){\line(1,0){2}}
\put(8,3){\line(2,-1){2}}\put(10,2){\line(3,-1){3}}
\put(1,2){\line(0,1){1}}\put(8,2){\line(0,1){1}}
\put(1,7){\line(0,1){1}}\put(1,7){\line(1,0){2}}
\put(1,7){\line(5,-1){5}}\put(1,8){\line(2,-1){2}}
\put(3,7){\line(3,-1){3}}\put(8,6){\line(0,1){2}}
\put(8,6){\line(4,1){4}}\put(8,8){\line(4,-1){4}}
\put(3.4,0.5){$\xi$}\put(5.7,0.5){$2\xi$}\put(0.6,1.9){$1$}\put(0.6,2.9){$2$}
\put(9.9,0.5){$\eta$}\put(12.5,0.5){$\xi+\eta$}\put(7.6,1.9){$1$}\put(7.6,2.9){$2$}
\put(11.9,5.5){$\xi$}\put(7.6,7.9){$2$}
\put(2.9,5.5){$\xi$}\put(5.5,5.5){$\xi+\eta$}\put(0.6,6.9){$1$}\put(0.6,7.9){$2$}
\put(3.7,5){(a)}\put(3.7,0){(c)}\put(10.7,5){(b)}\put(10.7,0){(d)}
\end{picture}
\caption{Tropical limits}\label{fig1}
\end{figure}

Let us prove statement (2) of Proposition
\ref{p6}. If $\bp\cap\widehat C^{(1)}\cap\widehat C^{(2)}=\emptyset$, then
the (immersed) curves $C^{(1)}=\bn(\widehat C^{(1)})$ and
$C^{(2)}=\bn(\widehat C^{(2)})$ intersect transversally and outside $\bz$, and the point
$\widehat z=\widehat C^{(1)}\cap\widehat C^{(2)}$ is mapped to a node of $C^{(1)}\cup C^{(2)}\setminus
\bz$. Then the uniqueness of the family $[\bn_t:\widehat C_t\to\Sigma,\bp_t]$,
$t\in(\C,0)$, and its smoothness follows from the standard properties of the deformation
smoothing out a node (see, for example, \cite[Lemma 11(ii)]{Itenberg_Kharlamov_Shustin:2012}).
Suppose now that the point $\widehat C^{(1)}\cap\widehat C^{(2)}$ belongs to $\bp$.
We prove statement (2) under condition $n=1$, leaving the case $n>1$ to the reader
as a routine generalization with a bit more complicated notations. Denote $\xi:=s^{(1)}_1=
=-D^{(1)}K_\Sigma$, $\eta:=s^{(2)}_1=-D^{(2)}K_\Sigma-1$.
We have three possibilities:
\begin{itemize}\item Suppose that $\xi<\eta$. In suitable coordinates
$x,y$ in a neighborhood of
$z_1=(0,0)$, we have
$$\alpha_1\equiv y-\lambda x^\eta\mod\fm_{z_1}^s,\quad C^{(1)}=\{y+x^\xi+\text{h.o.t.}=0\},
\quad C^{(2)}=\{y=0\}\ ,$$ where $\lambda\ne0$ is generic.
Without loss of generality, we can define
the family of arcs $(\bz_t,{\mathcal A}_t)_{t\in(\C,0)}$ by $\bz_t=(t,0)$,
${\mathcal A}_t=\{y\equiv\lambda(x-t)^\eta\mod\fm_{\bz_t}^s\}$ (cf. Lemma
\ref{le6}). The ideal $I_{z_1}$ from Lemma \ref{le6} can be expressed as
$\langle y^2,yx^{\xi-1},x^{\xi+\eta}\rangle$. Furthermore, by Lemma
\ref{le6}, for any family (\ref{ered6}), the curves
$C_\tau=\bn(\widehat C_\tau)\in|D|$ are given, in a neighborhood of $z_1$, by
\begin{eqnarray}&&y^2(1+O(x,y,\overline c))+yx^\xi(1+O(x,\overline c))+\sigma(\overline c)yx^{\xi-1}\nonumber\\
&&\qquad+\sum_{i=0}^{\xi-2}c_{i1}(\tau)yx^i
+\sum_{i=0}^{\xi+\eta-1}c_{0i}(\tau)x^i+O(x^{\xi+\eta},\overline c)=0\ ,\label{ered1}
\end{eqnarray}
where $\overline c$ denotes the collection of variables
$\{c_{i1},\ 0\le i\le \xi-2,\ c_{i0},\ 0\le i\le \xi+\eta-1\}$, the functions $c_{ij}(\tau)$ vanish
at zero for all $i,j$ in the summation range, and
$\sigma(0)=0$. Changing
coordinates $x=x'+t$, where $t=\varphi(\tau)$, we obtain the family of curves
\begin{eqnarray}&&y^2(1+O(x',y,t,\overline c))+y(x')^\xi(1+O(x',t,\overline c))+
\sigma' y(x')^{\xi-1}\nonumber\\
&&\qquad+\sum_{i=0}^{\xi-2}c'_{i1}y(x')^i
+\sum_{i=0}^{\xi+\eta-1}c'_{0i}(x')^i+t\cdot O((x')^{\xi+\eta},t,\overline c)=0\ ,\label{ered2}
\end{eqnarray} where
\begin{equation}\begin{cases}c'_{i1}&=\sum_{0\le u\le \xi-2-i}\binom{i+u}{i}t^uc_{i+u,1}
+\binom{\xi-1}{i}t^{\xi-1-i}\sigma\\
&\quad +t^{\xi-i}\left(\binom{\xi}{i}+O(t)\right)+O(t^{\xi-i},\overline c),\quad i=0,...,\xi-2,\\
c'_{i1}&=\sum_{u\ge0}\binom{i+u}{i}t^uc_{i+u,0},\quad i=0,...,\xi+\eta-1,\\
\sigma'&=\sigma+t(\xi+O(t,\overline c)).\end{cases}\label{ered4}\end{equation}
Next, we change coordinates $y=y'+\lambda(x')^\eta$ and impose the condition
\mbox{$(C_\tau\cdot(\bz_{\varphi(\tau)},{\mathcal A}_{\varphi(\tau)}))\ge \xi+\eta$}, which amounts in the
following relations
on the variables $\overline c'=\{c'_{i1},\ 0\le i\le \xi-2,\ c'_{i0},\ 0\le i\le \xi+\eta-1\}$:
\begin{equation}\begin{cases} &c'_{i0}=0,\ i=0,...,\eta-1,
\quad c'_{i0}+\lambda c'_{i-\eta,1}=0,\ i=l,...,\eta+\xi-2,\\
& c'_{\xi+\eta-1}+\lambda\sigma'=0\ .\end{cases}
\label{ered3}\end{equation}
The new equation for the considered family of curves is then
\begin{eqnarray}&F(x,y)=(y')^2(1+O(x',y',t,\overline c))+y'(x')^\xi(1+O(x',t,\overline c))\nonumber\\
&\qquad\qquad+(x')^{\xi+\eta}
(a+O(x',t,\overline c))+y'\left(\sum_{i=0}^{\xi-2}c'_{i1}(x')^i+\sigma'(x')^{\xi-1}\right)=0\ .
\label{ered8}\end{eqnarray}
with some constant $a\ne0$. Consider the tropical limit of the family (\ref{ered8})
(see \cite[Section 2.3]{Sh0} or Section \ref{sec2.1}). The
corresponding subdivision of $\Delta$ must be as shown in
Figure
\ref{fig1}(a). Indeed, first,
$c'_{01}\ne0$, since otherwise the curves $C_\tau$ would be singular at $\bz_t$ contrary to the
general choice of $(\bz_t,{\mathcal A}_t)$. Second, no interior point of the segment
$[(0,1),(\xi,1)]$ is a vertex of the subdivision, since otherwise the curves $C_\tau$ would have a positive
genus: the
tropicalization of $C_\tau$ would then be a tropical curve with a cycle which lifts to a handle of $C_\tau$
(cf. \cite[Sections 2.2 and 2.3, Lemma 2.1]{Sh0}). By a similar reason, the limit polynomial
$F^\delta_{\ini}/y'=\sum_{i=0}^\xi c^0_{i1}(x')^i$, where $\delta$ is the segment
$[(0,1),(\xi,1)]$, must be the $\xi$-th power of a
binomial.
The latter conclusion and relations (\ref{ered2}) and (\ref{ered4}) yield that $N_F(i,1)=\xi-i$
for $i=0,...,\xi$
and
$$c'_{i1}=t^{\xi-i}(c^0_{i1}+c''_{i1}(t)),\ i=0,...,\xi-2,\quad c'_{\xi+\eta-1,0}=t(c^0_{\xi+\eta-1,0}+
c''_{\xi+\eta-1,0})\ ,$$
where
$c^0_{i1}$, $i=0,...,\xi-2$, and $c^0_{\xi+\eta-1,0}$ are uniquely determined by the given data, the
functions
$c''_{i1}$, $0\le i\le \xi-2$, vanish at zero, and $c''_{\xi+\eta-1,0}$ is a function of
$t$ and $c''_{i1}$, $0\le i\le \xi-2$, that is determined by the given data and vanishes at zero too.
To meet the condition of
rationality of $C_\tau$ and to find the functions $c''_{i1}(t)$, $0\le i\le \xi-2$, we perform
the refinement
procedure
as described in \cite[Section 3.5]{Sh0}. It consists in further coordinate change and tropicalization,
in which one encounters a subdivision containing the triangle $\conv\{(0,0),(0,2),(\xi,1)\}$
(see Figure \ref{fig1}(b)).
The corresponding convex piece-wise linear function $N'$ is linear along that triangle and takes
values $N'(0,2)=N'(\xi,1)=0$, $N'(0,0)=\eta-\xi$. By \cite[Lemma 3.9 and Theorem 5]{Sh0}, there
are $\xi$ distinct solutions $\{c''_{i1}(t),\ 0\le i\le \xi-2\}$ of the rationality relation. More precisely,
the initial coefficient $(c''_{i1})^0$ is nonzero only for $0\le i\le \xi-2$, $i\equiv \xi\mod2$.
The common denominator of the values of $N'$ at these point is $\xi/d$, where $d=\text{gcd}(\xi,\eta)$, and
hence $c''_{i1}$
are analytic functions of $t^{d/\xi}$. It follows thereby that $t=\tau^{\xi/d}$.
\item Suppose that $\xi=\eta$. In this situation, the argument of the preceding case $\xi<\eta$ applies in
a similar way and, after the coordinate change $x=x'+t$, $y=y'+\lambda (x')^\xi$, leads to equation
(\ref{ered8}), whose Newton polygon is subdivided with a fragment $\conv\{(0,1),(0,2),(2\xi,0)\}$
on which the function $N_F$ is linear with values $N_F(0,2)=N_F(2\xi,0)=0$, $N_F(0,1)=\xi$.
By Lemma \ref{le7}, we get $\xi$ solutions $\{c'_{i1}(t),\ i=0,...,\xi-2\}$, which are analytic
functions of $t$.
Then, in particular, $t=\tau$.
\item Suppose that $\xi>\eta$. In suitable coordinates
$x,y$ in a neighborhood of
$z_1=(0,0)$, we have
$$\alpha_1\equiv y\mod\fm_{z_1}^s,\quad C^{(1)}=\{y+\lambda x^\xi+O(x^{\xi+1})=0\},
\quad C^{(2)}=\{y+x^\eta=0\}\ ,$$ where $\lambda\ne0$.
Without loss of generality, we can define
the family of arcs $(\bz_t,{\mathcal A}_t)_{t\in(\C,0)}$ by $\bz_t=(0,0)$,
${\mathcal A}_t=\{y\equiv tx^{\xi-1}\mod\fm_{z_1}^s\}$ (cf. Lemma
\ref{le6}). The ideal $I_{z_1}$ from Lemma \ref{le6} can be expressed as
$\langle y^2,yx^\xi,x^{\xi+\eta-1}\rangle$. Thus, by Lemma
\ref{le6}, for any family (\ref{ered6}), the curves
$C_\tau=\bn(\widehat C_\tau)\in|D|$ are given in a neighborhood of $z_1$ by
\begin{eqnarray}&&y^2(1+O(x,y,\overline c))+yx^\eta(1+O(x,\overline c))+
\lambda x^{\xi+\eta}(1+O(x,\overline c))\nonumber\\
&&\quad+\sigma(\overline c)x^{\xi+\eta-1}
+\sum_{i=0}^{\eta-1}c_{i1}(\tau)yx^i
+\sum_{i=0}^{\xi+\eta-2}c_{0i}(\tau)x^i=0\ ,\label{ered10}
\end{eqnarray}
where $\overline c$ now denotes the collection of variables
$\{c_{i1},\ 0\le i\le \eta-1,\ c_{i0},\ 0\le i\le \xi+\eta-2\}$, the functions $c_{ij}(\tau)$ vanish
at zero for all $i,j$ in the summation range, and
$\sigma(0)=0$. Inverting $t=\varphi(\tau)$, changing coordinates
$y=y'+tx^{\xi-1}$, and applying the condition
$(C_\tau\cdot{\mathcal A}_{\varphi(\tau)})\ge k+l$, we obtain an equation of the
curves $C_\tau$ in the form
\begin{eqnarray}&&F(x,y')=(y')^2(1+O(t,x,y',\overline c'))+y'x^\eta(1+O(t,x,\overline c'))\nonumber\\
&&\qquad\qquad
+\lambda x^{\xi+\eta}(1+O(t,x,\overline c'))+\sum_{i=0}^{\eta-1}c_{i1}(t)y'x^i
=0\ ,\label{ered11}
\end{eqnarray}
where $\overline c'=\{c_{i1},\ 0\le i\le \eta-1\}$, and the following relations must hold:
\begin{equation}\begin{cases} &c_{i0}=0,\ i=0,...,\eta-2,\\
&c_{i0}+tc_{i-\xi+1,1}=0,\ i=\xi-1,...,\eta+\xi-2,\\
& \sigma+t(1+O(t,\overline c'))=0\ .\end{cases}
\label{ered12}\end{equation}
By Lemma \ref{le4}(2), $\frac{\partial\sigma}{\partial c_{\eta-1,1}}(0)\ne0$. The rationality of
the curves $C_\tau$ yields that the subdivision $S_F$ of the Newton polygon
of $F(x,y')$ given by (\ref{ered11}) must contain two triangles
$\conv\{(0,1),(\eta,1),(0,2)\}$ and $\conv\{(0,1),(\eta,1),(\xi+\eta,0)\}$ (see Figure \ref{fig1}(d)),
and, furthermore, $F^\delta_{\ini}/y'$ must be the $\eta$-th power of a binomial, where $\delta
=[(0,1),(\eta,1)]$ (cf. the argument in the treatment of the case $\xi<\eta$ above). These two
conclusions and equations (\ref{ered12}) uniquely determine the initial coefficients
$c^0_{i1}$ as well as the values $N_F(i,1)=\eta-i$
for all $i=0,...,\eta-1$, and leave the final task to find the functions
$c''_{i1}(t)$, $i=0,...,\eta-2$, which appear in the expansion
$c_{i1}(t)=t^{\eta-i}(c^0_{i1}+c''_{i1}(t))$, $i=0,...,\eta-2$ (notice here that
the last equation in (\ref{ered12}) allows one to express $c''_{\eta-1,1}$
via $c''_{i1}$, $i=0,...,\eta-2$). To this extent, we again use the argument of the case
$\xi<\eta$, performing the refinement procedure along the edge $\delta=[(0,1),(\eta,1)]$
(see \cite[Section 3.5]{Sh0}) and apply the rationality requirement to draw the conclusion:
there are exactly $\eta$ families (\ref{ered6}), and, for each of them, $t=\tau^{\eta/d}$,
where $d=\text{gcd}\{\xi,\eta\}$.
\end{itemize}
Statement (2) of Proposition is proven. \proofend

\subsection{Families of curves and arcs on uninodal del Pezzo surfaces}
A smooth rational surface $\Sigma$ is called a uninodal del Pezzo surface if there exists a smooth rational
curve $E\subset\Sigma$ such that $E^2=-2$, and $-CK_\Sigma>0$ for each irreducible curve $C\subset \Sigma$
different from $E$. Observe that $EK_\Sigma=0$. Denote by $\Pic_+(\Sigma,E)\subset
\Pic(\Sigma)$ the semigroup generated by irreducible curves different from $E$. Assume that $\Sigma$
is of degree $1$ and fix $D\in\Pic_+(\Sigma,E)$ such that $-DK_\Sigma-1\ge3$.
Fix positive integers
$n\le-DK_\Sigma-1$ and \mbox{$s\gg-DK_\Sigma-1$}.

Accepting notations of Section \ref{sec2.2}, we introduce the set
$U^{im}(D,E)\subset \Arc^{\sm}_s(\ring\Sigma^n)$ is defined by the following conditions.
For any sequence $\bs=(s_1,...,s_n)\in\Z_{>0}^n$ summing up to $|\bs|\le s$ and for any
element $(\bz,{\mathcal A})\in U^{im}(D,E)$, where $\bz=(z_1,...,z_n)\in\ring\Sigma^n$,
$\bz\cap E=\emptyset$,
${\mathcal A}=(\alpha_1,...,\alpha_n)$, $\alpha_i\in
\Arc_s(\Sigma,z_i)$, the family
${\mathcal M}^{im}_{0,n}(\Sigma,D,\bs,\bz,{\mathcal A})$
is empty if $|\bs|\ge-DK_\Sigma$, and is finite
if $|\bs|=-DK_\Sigma-1$.
Furthermore, in the latter case,
all elements \mbox{$[\bn:\PP^1\to\Sigma,\bp]\in{\mathcal M}_{0,n}(\Sigma,D,\bs,\bz,{\mathcal A})$}
are represented by immersions $\bn:\PP^1\to\Sigma$ such that
$\bn^*(\alpha_i)=s_ip_i$, $1\le i\le n$, and $\bn^*(E)$ consists of $DE$ distinct points.

\begin{proposition}\label{p7}
The set $U^{im}(D,E)$ is Zariski open and dense in $\Arc^{\sm}_s(\ring\Sigma^n)$.
\end{proposition}

{\bf Proof.}
The statement that $U^{im}(D)$ is Zariski open and dense in $\Arc^{\sm}_s(\ring\Sigma^n)$
can be proved in the same way as Proposition \ref{p4}(1). We will show that
$U^{im}(D,E)$ is dense in $U^{im}(D)$, since the openness of $U^{im}(D,E)$ is evident.
For, it is enough to show that any immersion $\bn:\PP^1\to\Sigma$ such that
$\bn_*(\PP^1)=D$ can be deformed into an immersion with an image transversally crossing $E$
at $DE$ distinct points.

Suppose, first, that a generic element $[\bn:\PP^1\to\Sigma]\in{\mathcal M}_{0,0}(\Sigma,D)$
is such that the divisor $\bn^*(E)\subset\PP^1$ contains an $m$-multiple point, $m\ge2$.
Since $\dim{\mathcal M}_{0,0}(\Sigma,D)=-DK_\Sigma-1\ge3$, we fix the images of
$-DK_\Sigma-2$ points $p_i$, $i=1,...,-DK_\Sigma-2$, obtaining a one-dimensional subfamily
of ${\mathcal M}_{0,0}(\Sigma,D)$, for which one derives a contradiction by Lemma \ref{gus}(iii):
$$D^2\ge(D^2+DK_\Sigma+2)+(-DK_\Sigma-2)+(m-1)=D^2+m-1>D^2\ .$$
Hence, for a generic $[\bn:\PP^1\to\Sigma]\in{\mathcal M}_{0,0}(\Sigma,D)$, the divisor
$\bn^*(E)$ consists of $DE$ distinct points. Suppose that $m\ge2$ of them are mapped
to the same point in $E$. Fixing the position of that point on $E$, we define a subfamily
$V\subset{\mathcal M}_{0,0}(\Sigma,D)$ of dimension
$$\dim V\ge\dim{\mathcal M}_{0,0}(\Sigma,D)-1=-DK_\Sigma-2\ge2\ .$$
As above, we fix the images of $-DK_\Sigma-3$ additional point of $\PP^1$ and end up
with a contradiction due to Lemma \ref{gus}(ii):
$$D^2\ge(D^2+DK_\Sigma+2)+(-DK_\Sigma-3)+m=D^2+m-1>D^2\quad\text{\proofend}$$

Let ${\mathfrak X}\to(\C,0)$ be a smooth flat family of smooth rational
surfaces such that ${\mathfrak X}_0=\Sigma$ is a nodal del Pezzo surface with the $(-2)$-curve $E$,
and ${\mathfrak X}_t$, $t\ne0$, are del Pezzo surfaces.
We can naturally identify $\Pic({\mathfrak X}_t)\simeq\Pic(\Sigma)$,
$t\in(\C,0)$. Fix a divisor
$D\in\Pic_+(\Sigma,E)$ such that $-DK_\Sigma-1\ge3$. Given $n\ge1$ and $s\gg
-DK_\Sigma-1$, fix a vector $\bs\in\Z_{>0}^n$ such that $|\bs|=-DK_\Sigma-1$.
Denote by $\Arc^{\sm}_s({\mathfrak X})\to{\mathfrak X}\to(\C,0)$ the bundle with fibres $\Arc^{\sm}_s({\mathfrak X}_t)$, $t\in(\C,0)$.
Pick $n$ disjoint smooth sections $z_1,...,z_n:(\C,0)\to{\mathfrak X}$ covered by
$n$ sections $\alpha_1,...,\alpha_n:(\C,0)\to\Arc^{\sm}_s({\mathfrak X})$ such that
$(\bz(0),{\mathcal A}(0))\in U^{im}(\Sigma,E)$, and $(\bz(t),{\mathcal A}(t))\in
U^{im}({\mathfrak X}_t)$, $t\ne0$.

\begin{proposition}\label{p8}
Each element $[\nu:\widehat C\to\Sigma,\bp]
\in\overline{{\mathcal M}_{0,n}(\Sigma,D,\bs,\bz(0),{\mathcal A}(0))}$ such that
\begin{itemize}\item either $\widehat C\simeq\PP^1$, or $\widehat C=
\widehat C'\cup\widehat E_1\cup...\cup\widehat E_k$ for some $k\ge1$,
where $\widehat C'\simeq\widehat E_1\simeq...\simeq\widehat E_k\simeq\PP^1$,
$\widehat E_i\cap\widehat E_j=\emptyset$ for all $i\ne j$, and
$\#(\widehat C'\cap \widehat E_i)=1$ for all $i=1,...,k$;
\item $\bp\subset\widehat C'$ and $[\nu:\widehat C'\to\Sigma,\bp]\in{\mathcal M}^{im}_{0,n}(\Sigma,D-kE,\bs,\bz(0),{\mathcal A}(0))$, and
    each of the $\widehat E_i$ is isomorphically taken onto $E$;
\end{itemize} extends
to a smooth family $[\nu_t:\widehat C_t\to{\mathfrak X}_t,\bz(t)]\in \overline{{\mathcal M}_{0,n}}({\mathfrak X}_t,D,\bs,\bz(t),{\mathcal A}(t))$, $t\in(\C,0)$,
where $\widehat C_t\simeq\PP^1$ and $\nu_t$ is an immersion for
all $t\ne 0$, and, furthermore, each element of ${\mathcal M}_{0,n}
({\mathfrak X}_t,D,\bs,\bz(t),{\mathcal A}(t))$, $t\in(\C,0)\setminus\{0\}$ is included into some of the above families.
\end{proposition}

{\bf Proof.}
The statement follows from \cite[Theorem 4.2]{Vakil:2000} and from Proposition \ref{p7}, which applies to all
divisors $D-kE$, since $-(D-kE)K_\Sigma=-DK_\Sigma$ for any $k$.
\proofend

\section{Proof of Theorem \ref{t1}}

By blowing up
additional real points if necessary, we reduce the problem to consideration of
del Pezzo surfaces $X$ of degree $1$.

\smallskip

(1) To prove the first statement of Theorem \ref{t1} it is enough to consider
only del Pezzo surfaces satisfying property (GDP) introduced in Section \ref{secGDP} (cf. \cite[Lemma 17]{Itenberg_Kharlamov_Shustin:2012})
and real divisors satisfying $-DK_X-1\ge3$
(cf. Remark \ref{r1}(1)). So, let a real del Pezzo surface $X$ satisfy property (GDP)
and have a non-empty real part. Let $F\subset\R X$ be a connected component.
Denote by ${\mathcal P}_{r,m}(X,F)$ the set of sequences $(\bz,\bw)$ of $n=r+2m$ distinct points
in $\Sigma$ such that $\bz$ is a sequence of $r$ points
belonging to the component $F\subset\R X$, and $\bw$ is a sequence of
$m$ pairs of complex conjugate points. Fix an integer $s\gg-DK_X$ and
denote by $\R\Arc^{\sm}_s(X,F,r,m)\subset\Arc^{\sm}_s(\ring X^n)$ the space of
sequences of arcs $({\mathcal A},{\mathcal B})$ centered at
$(\bz,\bw)\in{\mathcal P}_{t,m}(X,F)$ such that ${\mathcal A}
=(\alpha_1,...,\alpha_r)$ is a sequence
of real arcs
$\alpha_i\in\Arc_s(X,z_i)$, $z_i\in\bz$, $i=1,...,r$, and
${\mathcal B}=(\beta_1,\overline\beta_1,...,\beta_m,\overline\beta_m)$ is a sequence of
$m$ pairs of complex conjugate arcs, where $\beta_i\in\Arc_s(X,w_i)$, $\overline\beta_i\in
\Arc_s(X,\overline w_i)$, $i=1,...,m$, and $\bw=(w_1,\overline w_1,...,w_m,\overline w_m)$.

We join two elements 
of $\R\Arc_s(X,F,r,m)\cap U^{im}(D)$ by a smooth real analytic path $\Pi
=\{(\bz_t,\bw_t),({\mathcal A}_t,{\mathcal B}_t)\}_{t\in[0,1]}$ in $\R\Arc_s(X,F,r,m)$
and show that along this path, the function $W(t):=W(X,D,F,\varphi,(\bk,\bl),(\bz_t,\bw_t),({\mathcal A}_t,
{\mathcal B}_t))$, $t\in[0,1]$, remains constant. By Propositions \ref{p4} and \ref{p6}, we need only to verify
the required constancy when the path $\Pi$ crosses sets $U^{im}_+(D)$, $U^{sing}_1(D)$, $U^{sing}_2(D)$,
$U^{mt}(D)$, and $U^{red}(D)$ at generic elements of their components of
codimension $1$ in $\Arc^{\sm}_s(\ring X^n)$. Let $t^*\in(0,1)$ correspond to the intersection of $\Pi$
with some of these walls.

If is clear that crossing of the wall $U^{sm}_+(D)\cap\R\Arc_s(X,F,r,m)$ does not affect
$W(X,D,F,\varphi,(\bk,\bl),(\bz_t,\bw_t),({\mathcal A}_t,
{\mathcal B})_t)$.

The constancy of $W(t)$ in a crossing of the wall $U^{sing}_1(D)\cap\R\Arc_s(X,F,r,m)$ follows from
Proposition \ref{p4}(3) and
\cite[Lemmas 13(2), 14 and 15]{Itenberg_Kharlamov_Shustin:2012}. The transversality hypothesis in
\cite[Lemma 15]{Itenberg_Kharlamov_Shustin:2012} can be proved precisely as
\cite[Lemma 13(1)]{Itenberg_Kharlamov_Shustin:2012}.

The constancy of $W(t)$ in a crossing of the wall $U^{sing}_2(D)\cap\R\Arc_s(X,F,r,m)$ follows from
Proposition \ref{p4}(4) and Lemma \ref{def}.

The constancy of $W(t)$ in a crossing of the wall $U^{mt}(D)\cap\R\Arc_s(X,F,r,m)$ follows from
Propositions \ref{p4}(5) and \ref{p5}. Indeed, by Proposition \ref{p5} exactly one real element of the set
${\mathcal M}_{0,n}(X,D,(\bk,\bl),(\bz_t,\bw_t),({\mathcal A}_t,{\mathcal B}_t))$ undergoes a bifurcation. Furthermore,
the ramification points of the degenerate map $\bn:\PP^1\to X$ are complex conjugate. Hence, the real part of a close
curve doubly covers the real part of $C=\bn(\PP^1)$, which means that the number of solitary nodes is always even.

At last, the constancy of $W(t)$ in a crossing of the wall $U^{red}(D)\cap\R\Arc_s(X,F,r,m)$ we derive from
Proposition \ref{p6}. Notice that the points $p_1\in\widehat C$ and $z_1\in X$ must be real, and hence
the cases (1ii) and (1iii) are not relevant, since we have the lower bound $-kD'K_X\ge2k\ge4$ contrary to
(\ref{e33}). In the case (1i) we use Proposition \ref{p6}(2): \begin{itemize}\item if $\bp\cap\widehat C^{(1)}\cap\widehat C^{(2)}=\emptyset$,
then the germ of the real part of the family (\ref{ered6}) is isomorphically mapped onto the germ $(\R,t^*)$ so that
the central curve deforms by smoothing out a node both for $t>t^*$ and $t<t^*$, and hence $W(t)$ remains
unchanged;
\item if $\bp\cap\widehat C^{(1)}\cap C^{(2)}=\{p_1\}$, then $p_1\in\PP^1$ and $z_1\in X$ must be real, and hence
$\xi+\eta$ must be odd, in particular,
$d=\gcd\{\xi,\eta\}$ is odd too, where $\xi=s^{(1)}_1$, $\eta=s^{(2)}_1$; if
$\kappa=\min\{\xi,\eta\}$ is odd, then
the real part of each real family (\ref{ered6}) is homeomorphically mapped onto
the germ $(\R,t^*)$, and, in the deformation
of the central curve both for $t>t^*$ and $t<t^*$, one obtains in a neighborhood
of $z_1$ an even number of real solitary nodes,
which follows from Lemma \ref{le7}(2); if $\kappa$ is even, then either the real
part
of a real family (\ref{ered6}) is empty, or the real part of a real family
(\ref{ered6}) doubly covers one of the
halves of the germ $(\R,t^*)$, so that in one component of $(\R,t^*)
\setminus\{t^*\}$ one has no real curves in the
family (\ref{ered6}), and in the other component of $(\R,t^*)\setminus\{t^*\}$
one has a couple or real curves,
one having an odd number $\kappa-1$ real solitary nodes, and the other having
no real solitary nodes (see Lemma \ref{le7}(2)),
and hence $W(t)$ remains constant in such a bifurcation.
\end{itemize}

\smallskip

(2) By \cite[Proposition 1]{Itenberg_Kharlamov_Shustin:2012}, in a generic one-dimensional
family of smooth rational surfaces of degree $1$ all but finitely many of them are del Pezzo and the
exceptional one are uninodal. Hence, to prove the second statement of Theorem \ref{t1} it is enough to
establish the constancy of $$W(t)=W({\mathfrak X}_t,D,F_t,\varphi,(\bk,\bl),(\bz(t),\bw(t)),
({\mathcal A}(t),{\mathcal B}(t)))$$
in germs of real families ${\mathfrak X}\to(\C,0)$ as in Proposition \ref{p8}, where the parameter is restricted to $(\R,0)\subset(\C,0)$. It follows from Proposition \ref{p8} that the number of
the real curves in count does not change, and real solitary nodes are not involved in the bifurcation.
Hence $W(t)$ remains constant.

\section{Examples}\label{sec2}

We illustrate Theorem \ref{t1} by a few elementary examples. Consider the case of plane cubics,
for which new invariants can easily be computed via integration with respect to the
Euler characteristic in the style of \cite[Proposition 4.7.3]{DK}.

Let $r_1+3r_3+2(m_1+2m_2+3m_3+4m_4)=8$, where $r_1,r_3,m_1,m_2,m_3,m_4\ge0$. Define integer vectors
$\bk=(r_1\times1,r_3\times 3)$, $\bl=(m_1\times1,m_2\times2,m_3\times3,m_4\times 4)$. Denote
by $L$ the class of line in $\Pic(\PP^2)$. Then
$$W(\PP^2,3L,(\bk,\bl))=r_1-r_3\ .$$
As compared with the case of usual Welschinger invariants, in the real pencil of plane cubics
meeting the intersection conditions with a given collection of arcs, in addition to real rational
cubics with a node outside the arc centers, one encounters rational cubics with
a node at the center of an arc of order $3$. Notice that this real node is not solitary since one of its local branches must be quadratically tangent to the given arc.
We also remark that, in a similar computation for a collection of arcs containing a real arc of order $2$,
one also encounters rational cubics with a node at the center of such an arc, but
this node can be solitary or non-solitary depending on the given collection of arcs, and hence the
count or real rational cubics will also depend on the choice of a collection of arcs.

Of course, the same argument provides formulas for invariants of any real del Pezzo surface and $D=-K$,
or, more generally, for each effective divisor with $p_a(D)=1$.

We plan to address the computational aspects in detail in a forthcoming paper.

\medskip

{\small
\section*{Acknowledgements} The author has been supported by the grant no. 1174-197.6/2011
from the German-Israeli Foundations, by the grant no.
176/15 from the Israeli Science Foundation, and by a grant from the Hermann-Minkowski-Minerva
Center for Geometry at the Tel Aviv University.}

\bigskip

{\it Address}: School of Mathematical Sciences, Tel Aviv
University, Ramat Aviv, 69978 Tel Aviv, Israel.
{\it E-mail}: {\tt shustin@post.tau.ac.il}

\end{document}